\documentclass[11pt, oneside]{amsart}
 \usepackage[text={5.58in,8.5in},centering,letterpaper,dvips]{geometry}
\usepackage{graphicx}
\usepackage{amsfonts}
\usepackage[dvipsnames]{xcolor}
\usepackage{subcaption}
\usepackage{epsf}
\usepackage{amssymb}
\usepackage{amsmath}
\usepackage{amscd}
\usepackage{mathtools}
\usepackage{pdfpages}
\usepackage{fancyhdr}
\usepackage{setspace}
\usepackage{soul} % for \ul \so \hl
\usepackage[all]{xy}
\usepackage{verbatim}
\usepackage{enumerate}
\usepackage[colorlinks=true, urlcolor=NavyBlue, linkcolor=NavyBlue, citecolor=NavyBlue]{hyperref}
    \usepackage{tikz}
\usetikzlibrary{  cd,  calc, positioning,  fit, arrows,  decorations.pathreplacing, decorations.markings,  shapes.geometric, backgrounds,bending}
\usepackage{tikzsymbols}

\theoremstyle{theorem}
\newtheorem{theorem}{Theorem}[section]
\newtheorem{proposition}[theorem]{Proposition}
\newtheorem{lemma}[theorem]{Lemma}
\newtheorem{question}[theorem]{Question}
\newtheorem{corollary}[theorem]{Corollary}

\makeatletter
\newtheorem*{rep@theorem}{\rep@title}
\newcommand{\newreptheorem}[2]{%
\newenvironment{rep#1}[1]{%
 \def\rep@title{#2 \ref{##1}}%
 \begin{rep@theorem}}%
 {\end{rep@theorem}}}
\makeatother

\newreptheorem{theorem}{Theorem}
\newreptheorem{lemma}{Lemma}
\newreptheorem{question}{Question}
\newreptheorem{corollary}{Corollary}
\newreptheorem{proposition}{Proposition}

\theoremstyle{definition}
\newtheorem{definition}[theorem]{Definition}
\newtheorem{remark}[theorem]{Remark}
\newtheorem{construction}[theorem]{Construction}

\newtheorem{example}[theorem]{Example}

\makeatletter
\def\@seccntformat#1{%
  \protect\textup{\protect\@secnumfont
    \ifnum\pdfstrcmp{subsection}{#1}=0 \bfseries\fi% subsection # in \bfseries
    \csname the#1\endcsname
    \protect\@secnumpunct
  }%
}  
\makeatother

\newcommand{\Z}{\mathbb{Z}}

\newcommand{\Spin}{\text{Spin}}

\begin{document}

\rhead{\thepage}
\lhead{\author}
\thispagestyle{empty}

%\tableofcontents
%\listoffigures

\raggedbottom
\pagenumbering{arabic}
\setcounter{section}{0}

%%%%%%%%%%%%%%%%%%%%%%%%%%%%%%%%%%%%%%%%%%%%%%%%%%%%%%%%
%%%%%%%%%%%%%%%%%%%%%%%%%%%%%%%%%%%%%%%%%%%%%%%%%%%%%%%%
%%%%%%%%%%%%%%%%%%%%%%%%%%%%%%%%%%%%%%%%%%%%%%%%%%%%%%%%

\title{Constructing knots with low rational genera}
%A new construction of knots bounding small genus surfaces in rational homology balls
%Kind of a mouthful....
% My working title was "A Menagerie of Rational and integral slice surfaces"

%"Knots of low rational genus" is a fine title
%\date{\today}

\author{Clayton McDonald}
\address{HUN-REN Alfréd Rényi Institute of Mathematics}
\email{claytkm@gmail.com}
\urladdr{https://sites.google.com/view/claytkm/}

\author{Allison N.  Miller}
\address{Department of Mathematics \& Statistics, Swarthmore College, Swarthmore, PA 19081}
\email{amille11@swarthmore.edu}
\urladdr{https://sites.google.com/view/anmiller/}

\begin{abstract}
	We give a flexible construction for knots in the 3-sphere that bound surfaces of unexpectedly low genus in punctured open books on 3-manifolds.  We use this construction to give the first examples of knots whose genus differs in different $\mathbb{Z}/2\mathbb{Z}$ homology balls. We also establish that every knot bounds a M{\"o}bius band in a rational homology ball, and that there are knots whose genus in $T^4$ and $B^4$ differ arbitrarily. 
\end{abstract}

\maketitle

\section{Introduction}

The smooth 4-dimensional Poincar{\'e} conjecture (S4PC) asks whether every smooth 4-manifold that is homotopy equivalent to $S^4$ must be diffeomorphic to $S^4$, and remains one of the outstanding open problems in low-dimensional topology. By work of Freedman, S4PC is equivalent to the conjecture that the topological 4-manifold $S^4$ has only one smooth structure. 
One of the many challenges in attempting to disprove S4PC is that most smooth 4-manifold invariants 
%that one might wish to use to obstruct a homotopy 4-sphere from being diffeomorphic to $S^4$ 
vanish in the absence of  second homology. 
However, as suggested by Freedman, Gompf, Morrison, and Walker in~\cite{FGMW_Man_and_machine}, a homotopy sphere $S$ might be distinguished from $S^4$ by the existence of a knot $K \subset S^3$ that bounds a smoothly embedded disk in $S^\circ \vcentcolon= S \smallsetminus \text{int }B^4$, yet does not bound a smoothly embedded disk in $B^4$. 
Recent work of Manolescu, Marengon, and Piccirillo~\cite{ManolescuMarengonPiccirillo24} provided progress in line with this strategy by showing that the answer to `does a given knot bound a null-homologous smoothly embedded  disk' can distinguish between smooth structures on certain indefinite 4-manifolds. 
However, one indication of the difficulty of fully executing the strategy of ~\cite{FGMW_Man_and_machine}  is that even the following is open. 
\begin{question}~\label{quest:sliceinzhb4}
Is there a knot $K$ that bounds a smoothly embedded disk in an integer homology ball $W$ with $\partial W=S^3$, and yet is not smoothly slice? 
\end{question}
% To answer this question in the affirmative, one needs both a method to construct knots that are slice in integer homology balls without obviously being slice in $B^4$, and a potential obstruction to sliceness in $B^4$ that does not automatically vanish on such knots. Of existing sliceness obstructions, the only ones that are not known to vanish on knots bounding disks in integer homology balls are Rasmussen's $s$-invariant~\cite{Rasmussen10s} and its variations. 
%Second, there is no known construction of potential candidate knots.
On the other hand, if one weakens `integer homology ball' to `rational homology ball' in Question~\ref{quest:sliceinzhb4}, the answer is well known:  every strongly negative amphichiral knot bounds a smoothly embedded disk in a rational homology ball $W$ with $|H_1(W)|=2$.~\cite{Kawauchi1980, Kawauchi09}. Work of Levine further shows that every strongly negative amphichiral (SNA) knot bounds a disk in the same rational homology ball~\cite{Levine23}, which one can verify is diffeomorphic to $\text{Spin}(\mathbb{RP}^3)^\circ$. 
%Recent work of Lidman-Piccirillo~\cite{LidmanPiccirillo25} shows that the figure-eight knot does not bound a disk in a different rational homology ball with the same integral homology.
%as $\text{Spin}(\mathbb{RP}^3)^\circ$. 
Notably, in all of these examples the rational homology ball in question has even order first homology. One might therefore ask the following.

\begin{question}~\label{quest:sliceinqhb4odd}
Is there a knot $K$ that bounds a smoothly embedded disk in a rational homology ball $W$ with $\partial W=S^3$ and $|H_1(W)|$ odd, and yet is not smoothly slice?
\end{question}

In this paper, we give a new construction that takes as input a 3-manifold $M$ and a knot $J \subset M$, and produces families of links in $S^3$ with bounded genus in certain open books on $M$. 
This construction is described in more detail in Section \ref{sec:construction}.
We combine this construction with Casson-Gordon signature computations to answer the `genus version' of Question~\ref{quest:sliceinqhb4odd}.
\begin{theorem}\label{thm:mainthm}
    There exist infinitely many knots that bound smoothly embedded genus one surfaces in  rational homology balls $W$ with  $ \partial W=S^3$ and $|H_1(W)|$ odd,  and yet do not bound smoothly embedded genus one surfaces in $B^4$. 
\end{theorem}

\begin{remark}
The manifold $W$ we use to prove Theorem~\ref{thm:mainthm} is $\text{Spin}(L(3,1))^\circ$, where we recall
that for a closed rational homology 3-sphere $Y$,
\[\Spin(Y):=(Y^\circ \times S^1) \cup_{S^2 \times S^1} (S^2 \times D^2)\]
is a rational homology 4-sphere with $H_1(\Spin(Y))\cong H_1(Y)$. 
%We can also think of this as the open book on $Y$ with trivial monodromy.
The knots we use are of the following form, where $J_0,J_1,J_2,J_3,J_4$ are certain connected sums of torus knots. 
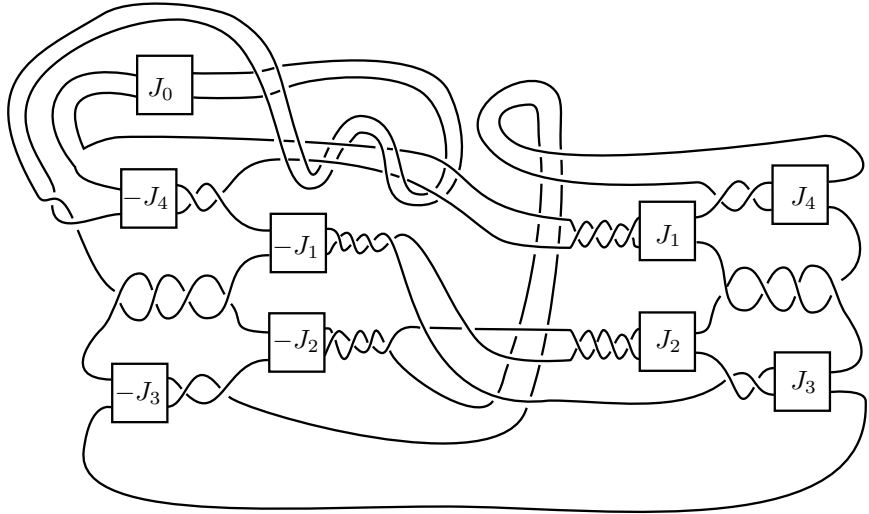
\begin{figure}[h!]
    \begin{tikzpicture}
    \node[anchor=south west,inner sep=0] at (0,0)  {\includegraphics[width=0.75\linewidth]{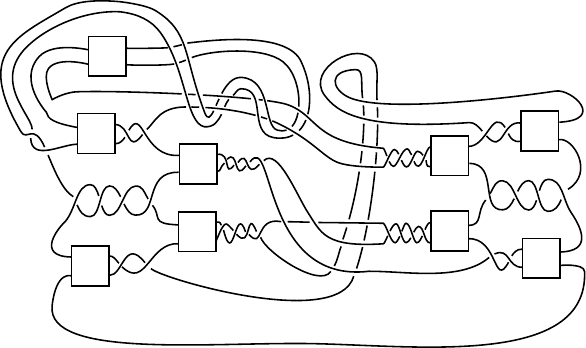}};
        \node at (1.75,1.6) {\small $-J_3$};
        \node at (1.85,4.2) {\small $-J_4$};
        \node at (2.05,5.65) {\small $J_0$};
        \node at (3.85, 2.3){\small $-J_2$};
        \node at (3.85, 3.55){\small $-J_1$};
        \node at (8.8, 2.3){\small $J_2$};
        \node at (8.8, 3.7){\small $J_1$};
        \node at (10.6, 1.75){\small $J_3$};
        \node at (10.6, 4.2){\small $J_4$};
    \end{tikzpicture}
    \caption{The knots $K(J_0,J_1,J_2,J_3,J_4)$ we use to prove Theorem~\ref{thm:mainthm}, subject to certain conditions on the companion knots $J_0,J_1,J_2,J_3,J_4$.}
    \label{fig:theexample}
\end{figure}

We will see that any knot $K=K(J_0,J_1,J_2,J_3,J_4)$ as in Figure~\ref{fig:theexample} has the property that $g_4(K \# K\#K) \leq 3$, and so most 4-genus bounds (including commonly used additive invariants like Rasmussen's $s$-invariant~\cite{Rasmussen10s}, the Heegaard-Floer $\tau$-invariant~\cite{OS03g4}, and the Levine-Tristram signature function~\cite{LevineSig,Tristram}) are unable to obstruct the 4-genus of $K$ from equaling one.
However, Casson-Gordon signatures, despite being additive with respect to connected sum of a (knot, character)-pair, have a sufficiently rich structure to establish Theorem~\ref{thm:mainthm}. Roughly, the fact that $g_4(K \#K \#K) \leq 3$ only forces certain sums of Casson-Gordon signatures of $K$ to be small, which can occur even if each individual signature is large. Similar arguments were used in~\cite{MillerAmphi} to show that there are strongly negative amphichiral (hence 2-torsion) knots with arbitrarily large 4-genera. 
\end{remark}

% While our methods do not provide an affirmative answer to Question~\ref{quest:sliceinqhb4odd}, we are able to address the following. 
% \begin{question}\label{quest:genusinqhb4odd}
% Is there a knot $K$ with smooth 4-genus $g$ that bounds a smoothly embedded surface of genus $h<g$ in a rational homology ball $W$ with $\partial W=S^3$ and $|H_1(W)|$ odd?
% \end{question}

In Section~\ref{sec:biggaps} we use connected sums of knots to obtain arbitrarily large gaps between the 4-genus and the minimal genus of a surface embedded in $\Spin(L(3,1))^\circ$.  
\begin{theorem}~\label{thm:biggaps}
For $m \in \mathbb{N}$, there is a knot $K_m$ with $g_4(K) \geq 3m$
that bounds a surface in $\Spin(L(3,1))^\circ$ of genus $2m$.
\end{theorem}
%Here, for a closed 4-manifold $Z$ and a knot $K$ in $S^3$, we let $g_Z(K)$ denote the minimal genus of a smoothly embedded surface in $Z\smallsetminus \text{int}(B^4)$ with boundary $K$. In this notation, $g_{S^4}$ is the classical 4-genus $g_4$.

\begin{remark}
    We expect that a direct analogue of our construction and arguments could be used to prove that for any $(p,q) \neq (1,0),(0,1)$, there is a knot $K$ that bounds a surface in $\Spin(L(p,q))^\circ$ with genus strictly smaller than $g_4(K)$. 
The chief difficulty in executing this lies in computing the Casson-Gordon signatures to obtain a lower bound on $g_4(K)$, since our construction will give increasingly complicated candidate knots as $p$ and $q$ grow. 
\end{remark}

We leverage a similar construction to obtain the following result. 
\begin{theorem}~\label{thm:mobiusinspinr3}
    Every knot in $S^3$ bounds a smoothly embedded M{\"o}bius band in $\text{Spin}(\mathbb{RP}^3)^\circ$.
\end{theorem}

% \begin{example}[Another perspective--Kirby calculus picture]
%     Everything bounds a mobius band in spin rp3.  In fact, also true for spin (L(4,1)), though with a slightly nonstandard diagram (and more). 
%     Want to have a collection of relators where every letter is used an even number of times. Genus is ``how many letters are used"/2. 
% \end{example}
Note that Theorem~\ref{thm:mobiusinspinr3} certainly does not hold in $B^4$~\cite{Yasuhara1996, MurakamiYasuhara, GilmerLivingston}. In fact, nonorientable genus can be arbitrarily large in $B^4$~\cite{Batson}.

Note that  a knot bounds a disk in $T^4$ if and only if it is smoothly slice in $B^4$; we prove a generalization of this in Proposition~\ref{prop:genuscollapset4}.
In contrast,  we use an alternate construction to show the following. 
\begin{theorem}~\label{thm:g4t4}
    The difference between the 4-genus of a knot and the minimal genus of a surface it bounds in $(T^4)^{\circ}$ can be arbitrarily large. 
\end{theorem}

The outline of our paper is as follows. In Section~\ref{sec:construction}, we develop our general construction, and go on to prove Theorem~\ref{thm:mobiusinspinr3} in Section~\ref{sec:spinrp3} and Theorem~\ref{thm:g4t4} in Section~\ref{sec:t4}. 
In Section~\ref{sec:CGstuff}, we perform the necessary Casson-Gordon signature computations to complete the proofs of Theorems~\ref{thm:mainthm} and~\ref{thm:biggaps}.  

\section*{Acknowledgements}
The authors would like to thank Marco Marengon for helpful comments on a draft of this article.
The first author would like to thank the Renyi Institute 2023 spring special semester "Singularities and Low Dimensional Topology", where most of the construction was conceived. The first author would also like to thank Maggie Miller for helpful conversations during the early stages of the construction.
The second author is also indebted to Anthony Conway for useful conversations about Casson-Gordon signatures and rational sliceness in summer 2019. 
The second author was partially supported by an AMS-Simons Research Enhancement Grant for PUI Faculty.

\section{Construction of surfaces}~\label{sec:construction}

In this section, we give a construction of knots in $S^3$ bounding surfaces of controlled genus in certain rational homology balls, and apply it to obtain various results. 

\subsection{General construction}
%In particular, we give a general family of rationally slice knots which includes the strongly negative amphichiral knots, as well as potential examples of knots whose genus in an integral ball is lower than its genus in $B^4$.

\subsubsection{Building surfaces in $\Spin(M)^\circ$}
We begin by describing how certain tangle cobordisms in $M^\circ \times I$ can be used to build surfaces in $\Spin(M)^\circ$. We first give the following definition.

\begin{definition}
 Let $T$ and $T'$ be tangles in $M^\circ$ with the same boundary $\partial T= \partial T'$.   
 A genus $g$ cobordism rel boundary between $T$ and $T'$ in $M^ \circ \times I$ is a properly embedded genus $g$ surface in $M^\circ \times I$ whose boundary consists of 
\[ -T \times\{0\} \,\cup\, \partial T \times [0,1] \,\cup\, T' \times \{1\}.\]
\end{definition}

The following proposition will be our key tool for building surfaces in $\Spin(M)^\circ$.
\begin{proposition}~\label{prop:keyconstruction}
Let $L\subset S^3$ be a link and $J \subset M$ be a knot. Let $T_J \subset M^\circ$ be the tangle obtained by removing a trivial arc ball from the pair $(M,J)$. 
If $L \# T_J$ is genus $g$ cobordant rel-boundary to $T_J$ in $M^\circ\times I$, then $L$ bounds a genus $g$ surface in $\Spin(M)^\circ$. 
%Namely, if $K\hash T_J$ is concordant to $T_J$ in $\mathring{M} \times I$, then $K$ is slice in punctured ${Spin(M)}$.

\end{proposition}

\begin{proof}
We decompose 
\begin{align*}
    \Spin(M) &= (M^\circ \times S^1) \, \cup_{S^2 \times S^1} \, (S^2 \times D^2),
\end{align*}
and further describe $M^\circ \times S^1= (M^\circ \times [-1,1])/\sim$, where as usual $(a,-1) \sim(a,1)$. 
Now choose $B^3 \subset M^\circ$ that is disjoint from $T_J$, in order to describe
\[ \Spin(M)^\circ= \left((M^\circ \times [-1,1]) \smallsetminus (\text{int}\, B^3 \times [-1/2,0])\right)/\sim \, \cup_{S^1 \times S^2}\, (S^2 \times D^2) \]
In Figure~\ref{fig:slicingprop}, we depict $M^\circ \times [-1,1]$ in black and the boundary of the removed $\text{int}(B^3) \times [-1/2,0]$ in red. We re-obtain $\Spin(M)^\circ$ by identifying $M^\circ \times \{-1\}$ with $M^\circ \times \{+1\}$ and then attaching $S^2 \times D^2$.

%This is to visualize the cobordism in the statement of the proposition inside of our manifold. We then describe the 4-dimensional puncture as $B^3 \times [-0.5,0]$, with $B^3$ as a submanifold of $\mathring{M}$ and $[-0.5,0]$ as a subinterval of $[-1,1]$ (see Figure \ref{fig:slicingprop} in red). We then imagine $K$ sitting in $B^3 \times \{0\}$. With this setup in place, we can use statements about cobordisms in $\mathring{M}\times I$ to prove statements about sliceness of knots in $S^3$ in $Spin(M)$.

Define
\[F_{[-1,0]}:= T_J \times [-1,0] \subset (M^\circ \times [-1,0]) \smallsetminus (\text{int}\,B^3 \times [-1/2,0]),\]
and let 
\[ F_{0}:= T_J \, \cup \, \mathfrak{b} \, \cup \, L\subset (M^\circ \smallsetminus \text{int}\,B^3) \times \{0\},\] where $\mathfrak{b}$ is a connect-sum band between $L$ and $T_J$. 
Finally, let $F_{[0,1]} \subset M^\circ \times [0,1]$ be the hypothesized cobordism from $K \# T_J$ to $T_J$. 

Then $F_{[-1,0]} \,\cup\, F_0 \,\cup \,F_{[0,1]}$ is a surface in $(M^\circ \times [-1,1]) \smallsetminus (\text{int}\,B^3\, \times [-1/2,0])$ with boundary  $T_J\times \{-1\} \,\cup\, \partial T_J \times [-1,1] \,\cup\, T_J\times \{1\}$ together with $L \times \{0\}$.

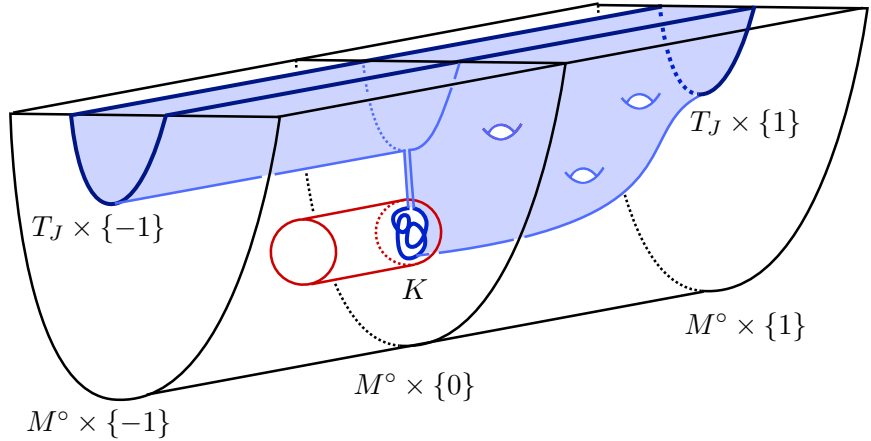
\begin{figure}[h]
    \centering
      \begin{tikzpicture}
        \node[anchor=south west,inner sep=0] at (0,0)  {\includegraphics[width=0.75\linewidth]{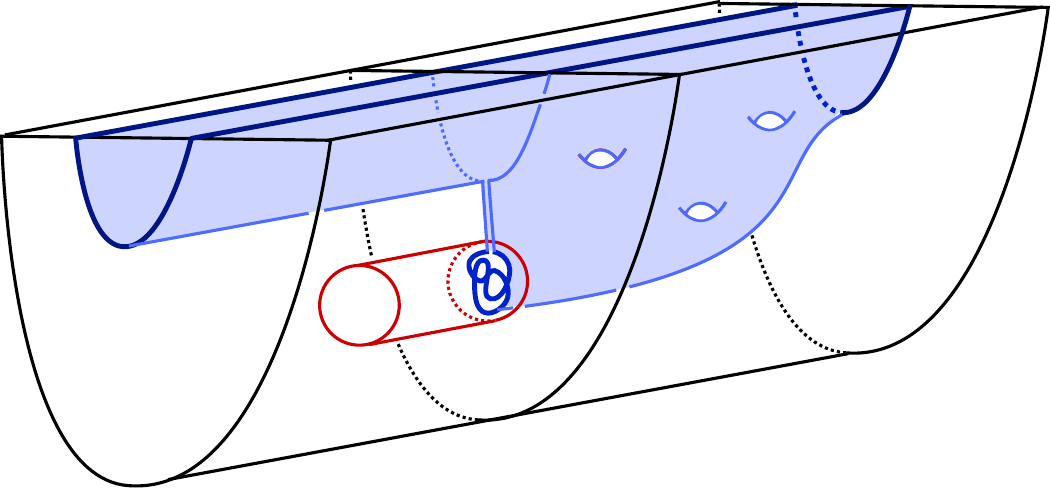}};
        \node at (1.2,-.3){$M^\circ \times \{-1\}$};
        \node at (1.2,2.3) {$T_J \times \{-1\}$};
        \node at (5.4, .2) {$M^\circ \times \{0\}$};
        \node at (9.8, 1){$M^\circ \times \{1\}$};
        \node at (9.8, 3.7){$T_J \times \{1\}$};
        %\draw[->] (9.8, 3.7)--(9.4,4.1);
        \node at (5.4, 1.5) {$K$};
        \end{tikzpicture}
    \caption{In black, we have $M^\circ\times [-1,1]$, with the extra 4-dimensional puncture in red as $B^3 \times [-1/2,0]$. In blue, we have the surface $F_{[-1,0]} \cup F_0 \cup F_{[0,1]}$.}
    \label{fig:slicingprop}
\end{figure}

Let $F$ denote the image of this surface in $(M^\circ \times S^1) \smallsetminus (\text{int}\,B^3 \times [-1/2,0])\subset \Spin(M)^\circ$. Note that $F$ is a genus $g$ surface with boundary given by  $\partial T_J \times S^1 \subset S^2 \times S^1$ and $L \subset \partial(B^3 \times [-1/2,0])$. Capping off the two components of $\partial T_J \times S^1$ with disks in $S^2 \times D^2$ gives our desired genus $g$ surface with boundary $K$ in $\Spin(M)^\circ$.
%We will now begin the proof by describing the portion of the surface in punctured $\mathring{M}\times I$. The surface, drawn schematically in blue in Figure \ref{fig:slicingprop}, is described as the following movie of tangles in 3-manifolds, starting from $t = 0$ for $t \in [-1,1]$
%We begin with $T_J \subset \mathring{M}$. We then see an extra puncture appear in the picture at $t = 0.5$, and at 0 the extra puncture disappears, but an extra link component of $K$ appears. we then band the two components together to get $T_J \# K$, and then attach the cobordism to $T_J$. Because the start and endpoints of this are both $T_J \subset \mathring{M}$, we can glue this up to form a genus g surface with three boundary components in $\mathring{M} \times S^1$. The three boundary components of this surface are $K$ and $\partial (T_J)\times S^1 \subset S^2\times S^1$. These $S^1$ factors bound discs in $S^2 \times D^2$, so we have the required genus g surface with boundary $K$.
\end{proof}

\begin{remark}
    One does not have to require a rel-boundary cobordism in the hypotheses of Proposition~\ref{prop:keyconstruction}.
    Given a surface in $M^\circ \times I$ with boundary $T \times \{0\} \, \cup\, \sigma \,\cup\, T' \times \{1\} $, where $\sigma$ is some 2-braid in $S^2 \times I$ connecting $\partial T \times \{0\}$ to $\partial T' \times \{1\}$, 
    we only need the induced braid on the boundary $S^1 \times S^2$ to be able to be capped off with discs.
    Happily, there are only two 2-braids in $S^1 \times S^2$, determined by their number of components, so to obtain our desired conclusion we need only check that $\alpha$ sends each endpoint of $T$ to same endpoint of $T'$, rather than swapping them. 
\end{remark}

\begin{remark}
    There are many potentially useful modifications of the above construction. For example, one could build surfaces in various open books by gluing the two ends of ${M}^\circ \times I$ according to some nontrivial mapping class $\phi$, which would then require $K \# T_J$ to be genus $g$ cobordant to $\phi(T_J)$. One could also replace $M^\circ \times I$ with some nontrivial self-cobordism of $M^\circ$. 
  %  A genus $g$ surface in ${M} \times I$ gives a corresponding genus $g+1$ surface in $M \times S^1$, or one can replace ${M} \times I$ with some nontrivial self-cobordism of $M$.
    %One must simply understand how to control the various extra added parameters to produce something coherent and useful to a new problem.
\end{remark}

%We will colloquially say that knots $K \subset S^3$ satisfying the conditions of Proposition~\ref{prop:keyconstruction} are {\bf dissolvable}, and call $J \subset M$ the {\bf solvent}.

It is not immediately clear that a surface satisfying the hypotheses of Proposition~\ref{prop:keyconstruction} occurs in a novel way. 
In particular, while the conclusion of Proposition~\ref{prop:keyconstruction} would also follow if we knew that $L$ bounded a genus $g$ surface in $(M^\circ\times I)^\circ$, we see that the presence of $T_J$ is necessary to obtain interesting examples as follows.  By taking the universal cover, one can see that $L$ bounds a  genus $g$ surface in $(M^\circ \times I)^\circ$ if and only if $L$ bounds a genus $g$ surface in $B^4$. 
That is, to effectively apply Proposition~\ref{prop:keyconstruction}, we need to  build tangle cobordisms from $T_J \#L$ to $L$ in $M^\circ \times I$ that do not arise from surfaces bounded by $L$ in $B^4$.

%As mentioned in the proof, if we instead require our local knot $K$ to simply be null cobordant in $M \times I$, then by taking the universal cover we see that this does not create any examples not already found in $S^3$ (see \cite{} for a full proof).
\subsubsection{Building surfaces in $M^\circ \times I$}

Traditionally, link cobordisms in $S^3 \times I$ are drawn as movies of links in $S^3$, with one of the following moves bridging the gap between frames:
\begin{enumerate}
\item isotopy
\item 0-handle attachment: the appearance of an split unknot.
\item 1-handle attachment: a band move.
\item 2-handle attachment: the disappearance of an split unknot.
\end{enumerate}

Modifying this approach to diagrammatically represent a tangle cobordism in $M \times I$   depends on a choice of diagrammatic representation for the 3-manifold $M$.
%For example, in a Dehn surgery diagram for a 3-manifold, the key difference is that isotopy can potentially involve sliding our knots and links over the surgery curves. Alternatively, we can phrase this as one extra move consisting of adding 0-handles to our diagram parallel to the meridians of our surgery curves.
We give a description in terms of Heegaard diagrams, recalling some definitions in order to fix notation. 
\begin{definition}
A Heegaard diagram is a triple $(\Sigma,\vec{\alpha} = \{\alpha_1, \dots \alpha_g\}, \vec{\beta} = \{\beta_1, \dots \beta_g\})$ of a closed genus $g$ surface $\Sigma$ along with two collections of disjoint simple closed curves $\vec{\alpha}$ and $\vec{\beta}$ on $\Sigma$ that each generate a $\Z^g$ summand of $H_1(\Sigma; \Z) \cong \Z^{2g}$. 
Using this triple, we can construct a closed 3-manifold in the following way.  First, take $\Sigma \times [-1,1]$ and attach thickened discs to $\Sigma \times \{-1\}$ along $\vec{\alpha}$, and to $\Sigma \times \{1\}$ along $\vec{\beta}$. The resulting 3-manifold will have two $S^2$ boundary components, which can be uniquely filled in with 3-balls.  
If the result is orientation preservingly homeomorphic to $M$, then we say $(\Sigma,\vec{\alpha}, \vec{\beta})$ is a \emph{Heegaard diagam for $M$}. We will sometimes use the decomposition $M= H_{\alpha} \cup \left(\Sigma \times [-1,1]\right) \cup H_{\beta}$, where $H_\alpha$ and $H_\beta$ correspond to the handlebodies described above. 
 \end{definition}

In our setting, the key fact we will use is that the $\alpha$-curves on $\Sigma \times \{-1\}$ and the $\beta$-curves on $\Sigma \times \{1\}$ all bound disks in $M$.
We begin with an example.

\begin{example}~\label{exl:fig8!}
    Let $M=\mathbb{RP}^3$, with the standard genus one Heegaard splitting. 
    We choose a basis for $H_1(T^2)$ such that the $\alpha$-curve represents $(2,1)$,  the $\beta$-curve represents $(0,1)$, and $J$ is the $(1,0)$ curve. 
    Recalling that we need to build a tangle cobordism in $M^\circ \times I$, we let $B^2$ be a small disk in $T^2$ that is disjoint from the $\alpha$- and $\beta$-curves. 
    Then our model for $M^\circ$ is
    \[M^\circ= H_\alpha \cup \left(\left( T^2 \times [-1,1]\right) \setminus \left(B^2 \times [-1/2,1/2]\right)\right) \cup H_\beta.
    \]
We think of $T^2$ as usual as the quotient of a square, in order to work diagrammatically, and think of the $\alpha$-curve as living in $T^2 \times \{1\}$, the arc $T_J$ in $T^2 \times \{0\}$, and the $\beta$-curve in $T^2 \times \{-1\}$. 

Our cobordism will be built starting with $T_J \times I$, adding four disks, and attaching four bands, as follows. At the beginning, as illustrated on the left of Figure~\ref{fig:fig8} we see only $T_J$. 
\begin{figure}[h!]
    \includegraphics[height=4cm]{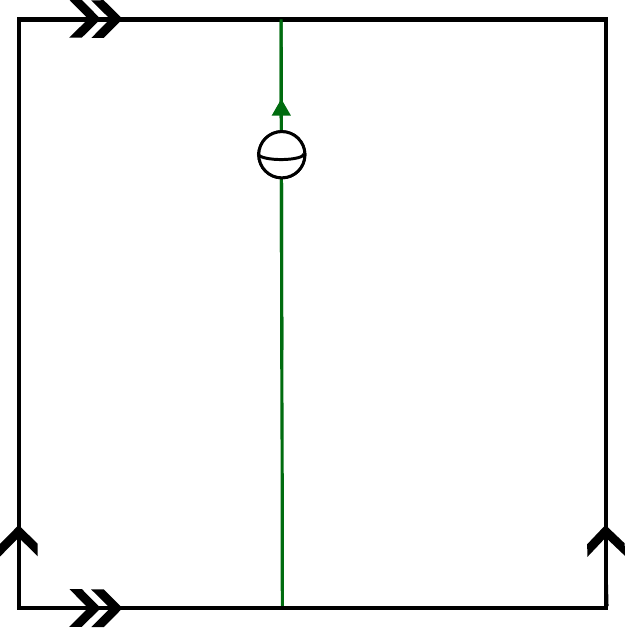}
    \quad
        \includegraphics[height=4cm]{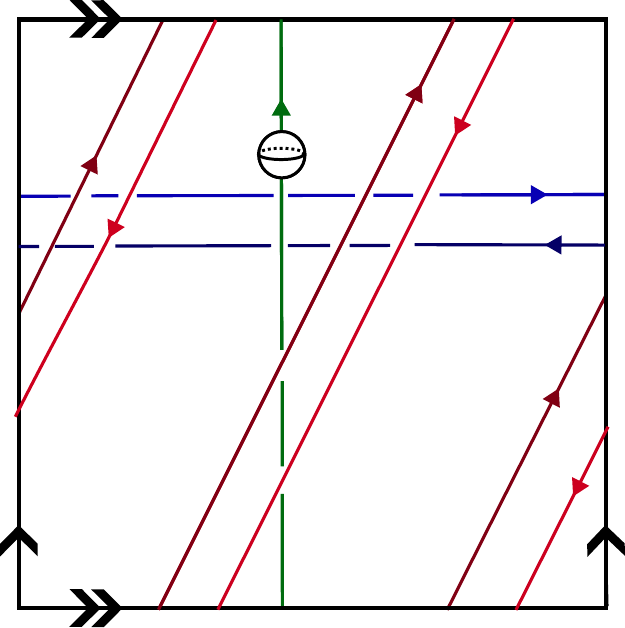}
        \quad
        \includegraphics[height=4cm]{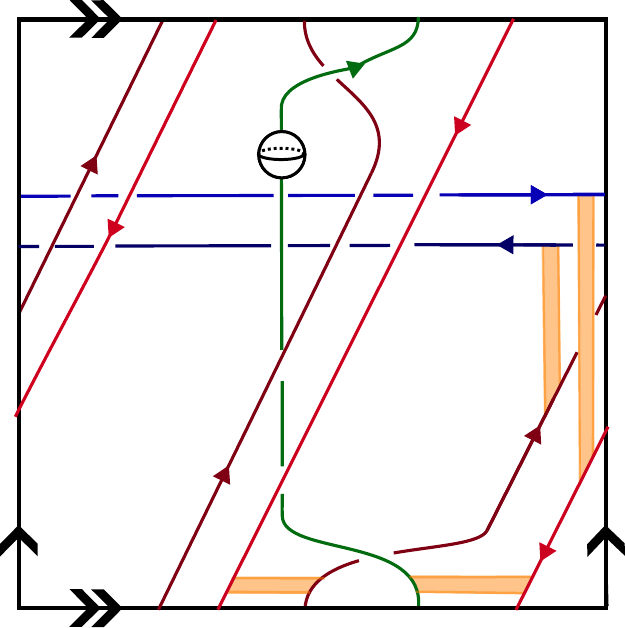}
    \caption{Left: The initial tangle $T_J$. Center: Four curves  appear. Right: Canceling half twists appear, together with four orange bands that will be added in the next step.}
    \label{fig:fig8}
\end{figure}
In the center of Figure~\ref{fig:fig8}, we see two oppositely oriented parallel copies of the $\alpha$-curves and two oppositely oriented parallel copies of the $\beta$-curves appear, corresponding to disks in the surface we are building in $M^\circ \times I$. The right of Figure~\ref{fig:fig8} is obtained by an isotopy in $M^\circ$.

We now attach the four orange bands shown on the right of Figure~\ref{fig:fig8}, and perform a small isotopy across the boundary of the fundamental domain $D$ to obtain the tangle shown on the left of Figure~\ref{fig:fig82}. After a further simplifying isotopy, we have the tangle shown in the center of Figure~\ref{fig:fig82}. 
    \begin{figure}[h!]
    \includegraphics[height=4cm]{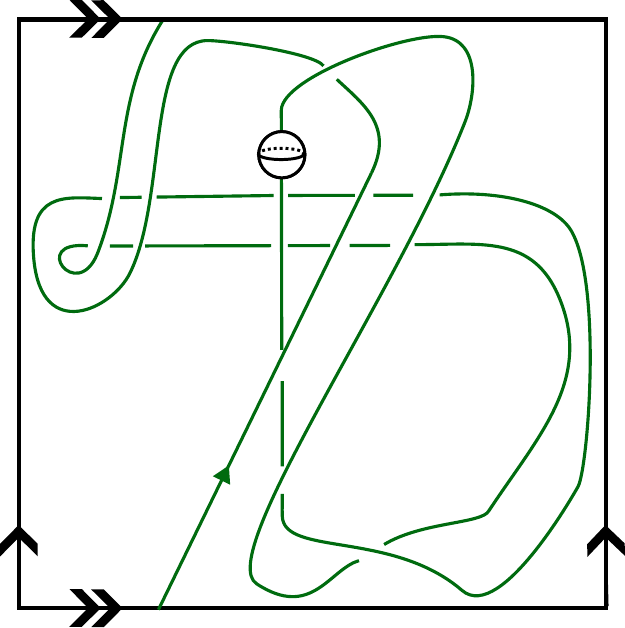}
    \quad
        \includegraphics[height=4cm]{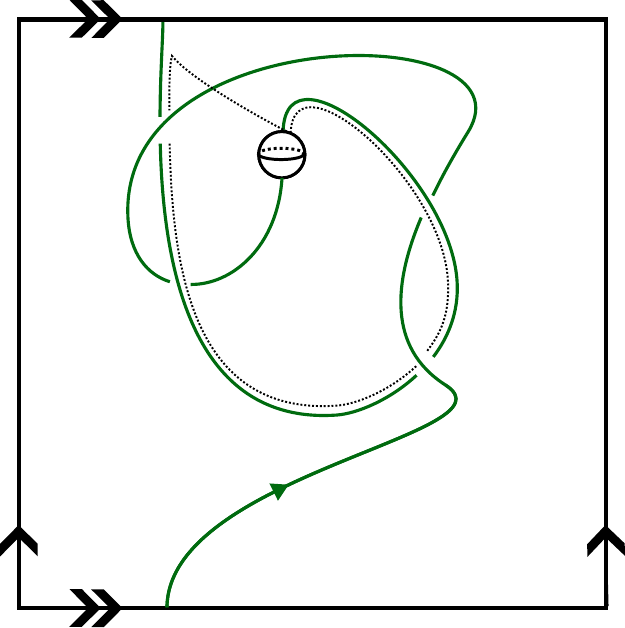}
        \quad
        \includegraphics[height=4cm]{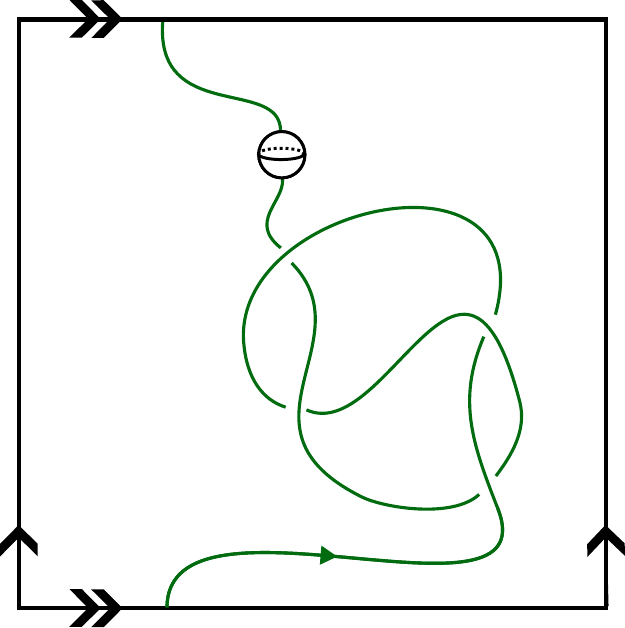}
    \caption{The tangle obtained by attaching bands (left), after isotopy (center), and after yet another isotopy (right).}
    \label{fig:fig82}
\end{figure}
We now claim that a further isotopy takes the center tangle to the one shown on the right, which is evidently $T_J \,\#\,K$ for $K$ the figure-eight knot. One way to see this is to imagine isotoping the removed $B^3$ in $M^\circ$ along the dotted loop, dragging the tangle with it for the first portion and then pushing it ahead along the final short `southeast' arc. This gives a self-homeomorphism of $(M, B^3)$, which gives a self-homeomorphism of $M^\circ$. Since the dotted loop is contained within a ball in $M$, this homeomorphism is isotopic to the identity. 
\end{example}

We note that it is not a coincidence that we obtain the figure eight knot here: we will shortly give a recipe for building a genus zero cobordism from $T_J$ to $T_J \, \#\, K$ for any strongly negative amphichiral knot $K$. However, first we generalize our cobordism building strategy for later use. 

\begin{construction}~\label{construction:surfacediagram}
Using a Heegaard diagram $(\Sigma, \vec{\alpha}, \vec{\beta})$ for $M$, we will choose $T_J \subset M^\circ$ and diagrammatically construct $L$ satisfying the conditions of Proposition~\ref{prop:keyconstruction} as follows. 

Let $B^2$ be a small open disk in $\Sigma$ that is disjoint from all $\alpha$- and $\beta$-curves. 
Our diagrams will live in $(\Sigma \times [-1,1]) \smallsetminus (B^2 \times [-1/2,1/2])\subset M^\circ$.  
For convenience, let $\Sigma^\circ= \Sigma \setminus B^2$. 
%The moves, as before, will be exactly the same as the normal 2-dimensional ambient handle attachments, but with an augmented set of 0- and 2-handles.

We begin by considering $\Sigma \times \{0\}$ as a quotient of a polygon $D$.
Choose a curve $J \subset \Sigma \times \{0\}$ with the properties that
\begin{enumerate}
    \item $J$ intersects the image of $\partial D$ in $\Sigma \times \{0\}$ in a single point that is not the image of a vertex of $D$, and 
    \item $J$ intersects $B^2 \times \{0\} \subset \Sigma \times \{0\}$ in a single arc. 
\end{enumerate}

Let 
\[T_J=(J \cap \Sigma^\circ) \times \{0\} \subset \Sigma^\circ \times \{0\} \subset M^\circ.\]
Now, add some number of parallel copies of the $\alpha$-curves (at the $\Sigma \times \{1\}$ level) and the $\beta$-curves (at the $\Sigma \times \{0\}$ level) to the diagram. This corresponds to introducing $n_\alpha+ n_\beta$ disks to the surface we're building in $M^\circ \times I$.

We now attach bands between the copies of the $\alpha$- and $\beta$-curves and $T_J$
to obtain a tangle $T_L \subset (\Sigma \times [-1,1]) \smallsetminus (B^2 \times [-1/2,1/2]) $ so the following hold:
\begin{enumerate}
    \item The intersection of $T_L$ with the image of $\partial D \times [-1,1]$ in $\Sigma \times [-1,1]$ is the same as that of 
    $T_J$ with the image of $\partial D \times \{0\}$. 
    \item The component of $T_L$ that intersects the image of $\partial D \times [-1,1]$ is the same component that meets $\partial B^2 \times [-1/2,1/2]$.
\end{enumerate}

% so as to remove geometric intersections with $\partial D$ and obtain a link $L$, so that exactly one intersection with $\partial D$ remains, and on the same component of $\partial D$ as $J \cap \partial D$.
% You must also choose bands such that the punctured component is the same component as the one that intersects $\partial D$.

Because of our geometric intersection constraints, we have that $T_L$ must be isotopic to $T_J \# L_J$ for some local link $L_J$, as the remaining knot is contained entirely within a 3-ball in $M^\circ$. 
Moreover, by construction the link $L_J$ has the property that $T_J$ and $T_J \, \# \, L_J$ cobound a surface in $M^\circ \times I$, whose genus can be computed by keeping track of the number of $\alpha$-curves, $\beta$-curves, and bands used. 
\end{construction}

\subsection{Examples in $\text{Spin}(\mathbb{RP}^3)$}
~\label{sec:spinrp3}

We now consider the local knots that arise from Construction~\ref{construction:surfacediagram} when $M= \Spin\left( \mathbb{RP}^3 \right)^\circ$.

    \begin{construction}~\label{cons:amphi}
    Start with the standard genus one Heegaard diagram of $\mathbb{RP}^3$, with a basis of $H_1(T^2)$ such that the $\alpha$- and $\beta$-curves represent the $(2,1)$ and $(0,1)$ curves respectively.
    Our base knot $J$ in the 3-manifold will be the $(1,0)$ curve. The left of Figure~\ref{fig:rp3} illustrates the tangle $T_J$ (green), together with $2k$ parallel copies of the $\alpha$- (red) and $2k$ parallel copies of the $\beta$-curves (blue), coming in oppositely oriented pairs. 
    We obtain the center and right of Figure~\ref{fig:rp3} by  isotopy of the $\alpha$ curves within $(\Sigma \smallsetminus \text{int}\, B^2) \times [.5,1]$. Here $\sigma$ is an arbitrary braid on $4k+1$ strands, and $\sigma^{-1}$ denotes its inverse in the braid group. 
    \begin{figure}[h!]
    \centering
% \begin{tikzpicture}
%         \node[anchor=south west,inner sep=0] at (0,0)  {\includegraphics[width=0.27\linewidth]{amphi1.pdf}};\end{tikzpicture}
  \begin{tikzpicture}
        \node[anchor=south west,inner sep=0] at (0,0)  {\includegraphics[width=0.3\linewidth]{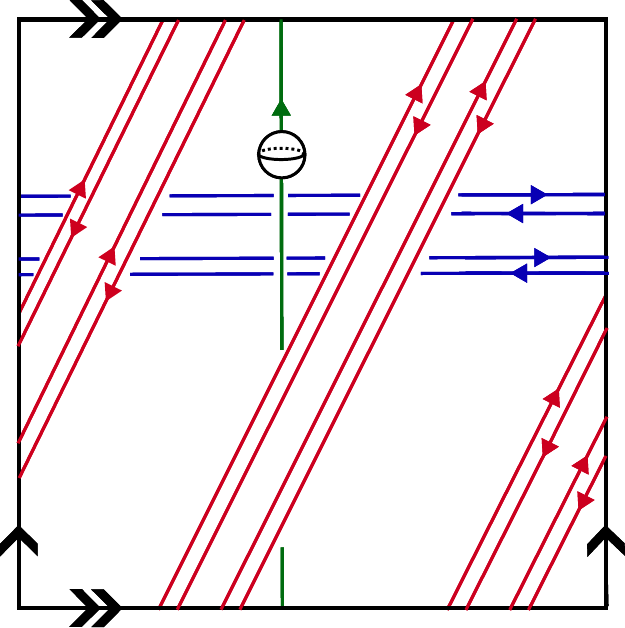}};\end{tikzpicture}
        \quad     
 \begin{tikzpicture}
 \node[anchor=south west,inner sep=0] at (0,0)  {\includegraphics[width=0.3\linewidth]{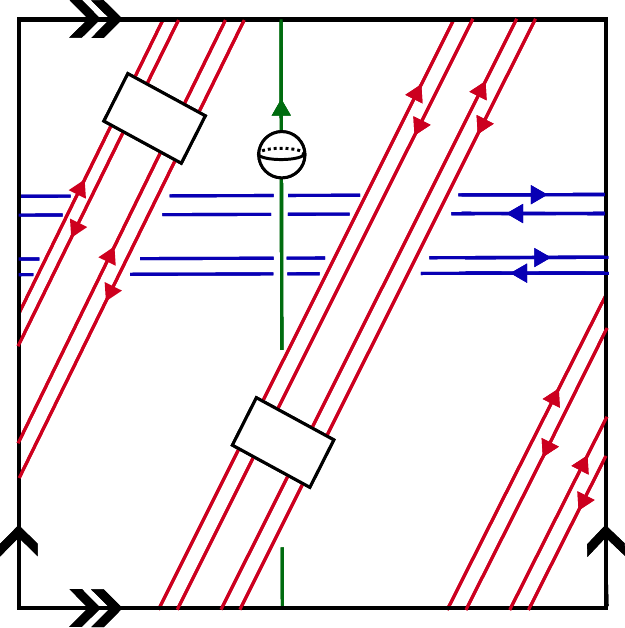}};
\node at (1.1,3.7)[rotate=-30]{ \tiny $-\tau$};
 \node at (2.1,1.35)[rotate=-30]{ \tiny$\tau$}; \end{tikzpicture}\quad
\begin{tikzpicture}
\node[anchor=south west,inner sep=0] at (0,0)  {\includegraphics[width=0.3\linewidth]{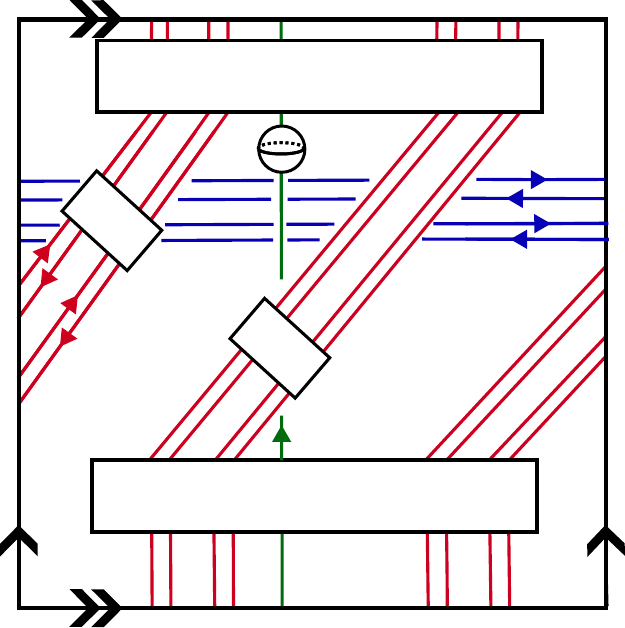}};
\node at (0.82,2.92)[rotate=-35]{ \tiny $-\tau$};
 \node at (2.1,2)[rotate=-35]{ \tiny$\tau$}; 
 \node at (2.3,4.05){$\sigma^{-1}$};
 \node at (2.2,.9){$\sigma$};
\end{tikzpicture}
\caption{(Left) In a neighborhood of the Heegaard torus, we see $J$ together with $k$ oppositely oriented pairs of $\alpha$- and $\beta$-curves. (Center) We have introduced a canceling pair of half full twists via isotopy. (Right) We have further introduced a braid $\sigma$ and its inverse $\sigma^{-1}$.}
\label{fig:rp3}
\end{figure}
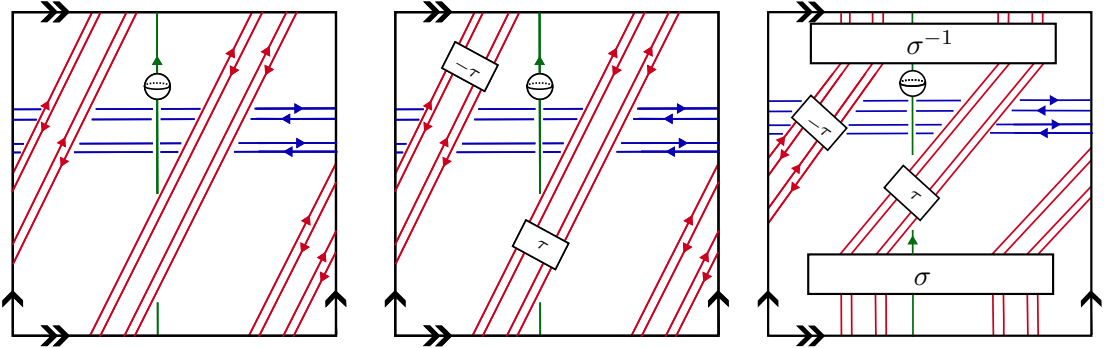

We would now like to attach bands to cancel as many intersection points of our curves with the boundary of our fundamental domain $D$ as possible. As shown in Figure~\ref{fig:fig8} for the case of $k=1$, we do so by attaching $2k$ bands, parallel to the right side of $D$, joining each $\alpha$-curve to a $\beta$-curve, and then attaching $2k$ more horizontal bands, parallel to the top side of $D$, joining 
$\alpha$-curves to each other (or, in one case, to $T_J$).
After a slight isotopy,  the result of banding removes almost all the intersection points of our curve with the boundary of $D$.

Note that after banding, as illustrated on the left of Figure~\ref{fig:rp32},  we have one component if and only if all of these band moves are orientation preserving. We will restrict to the case when this is true. We remark that this is why we chose $k$ pairs of oppositely oriented copies of the $\alpha$- and $\beta$-curves--otherwise, we would not be able to remove all but one intersection point via our banding.

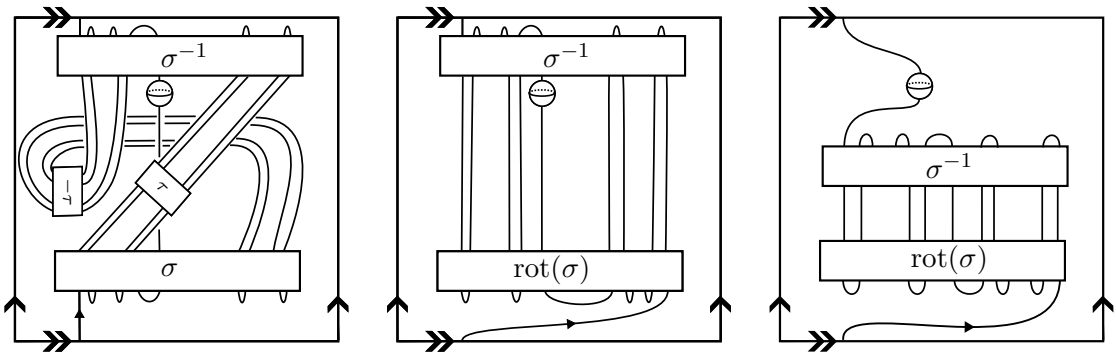
\begin{figure}[h!]
\begin{tikzpicture}
        \node[anchor=south west,inner sep=0] at (0,0)  {\includegraphics[width=0.3\linewidth]{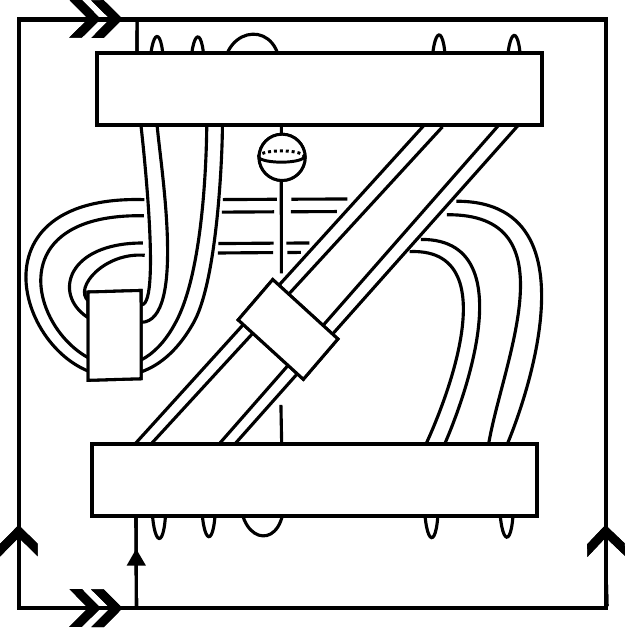}};
\node at (0.85,2.13)[rotate=-90]{ \tiny $-\tau$};
 \node at (2.1,2.2)[rotate=-45]{ \tiny$\tau$}; 
 \node at (2.2,1.1){$\sigma$};
 \node at (2.4,3.95){$\sigma^{-1}$};
        \end{tikzpicture}
        \quad
\begin{tikzpicture}
        \node[anchor=south west,inner sep=0] at (0,0)  {\includegraphics[width=0.3\linewidth]{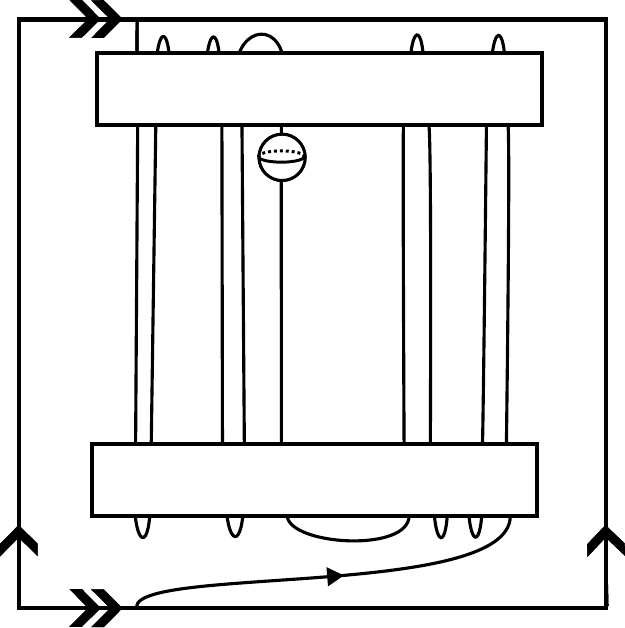}}; 
         \node at (2.2,1.1){$\text{rot}(\sigma)$};
 \node at (2.4,3.95){$\sigma^{-1}$};
 \end{tikzpicture}
        \quad
\begin{tikzpicture}
        \node[anchor=south west,inner sep=0] at (0,0)  {\includegraphics[width=0.3\linewidth]{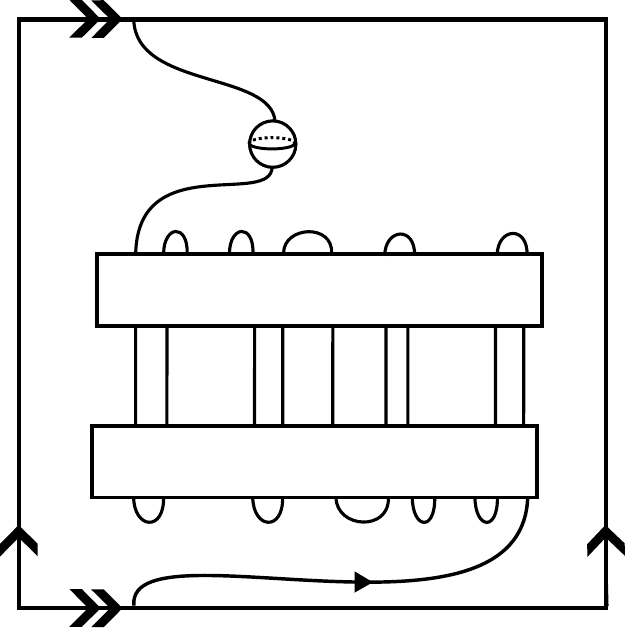}};
         \node at (2.4,1.23){$\text{rot}(\sigma)$};
 \node at (2.4,2.48){$\sigma^{-1}$};
 \end{tikzpicture}
\caption{We obtain the figure on the left from the rightmost part of Figure~\ref{fig:rp3} by attaching bands. We then obtain the center figure by an isotopy, and the right figure by a ``point pushing'' homeomorphism.} 
    \label{fig:rp32}
\end{figure}

By swinging the bottom part of the diagram around, we obtain the middle of
Figure~\ref{fig:rp32}, where we let $\text{rot}(\sigma)$ denote the braid obtained by rotating $\sigma$ by 180 degrees. 
We would now like to recognize the tangle $T$ shown on the center of Figure~\ref{fig:rp32} as $T_J \#K$ for a knot $K$. We will do this by recognizing the center  of Figure~\ref{fig:rp32} as isotopic to the right.  Without the puncture, this would be immediate. With the puncture, we follow an argument that is essentially identical to that used in Example~\ref{exl:fig8!}. That is,  we imagine isotoping the removed $B^3$ in $M^{\circ}$ along $T$ until it is just below the depicted intersection of $T$ with the top side of $D$, and then pushing it `southeast' to its original starting point. This gives a self-homeomorphism of $(M, B^3)$, which gives a self-homeomorphism of $M^\circ$ that sends $T$ to $T_J \# K$. Since the dotted loop is contained within a ball in $M$, this homeomorphism is isotopic to the identity.

If we look at the local knot $K$ we get from this procedure, illustrated in Figure~\ref{fig:rp33}, we can see it is strongly negative amphichiral.
    \end{construction}

We now combine Construction~\ref{cons:amphi} with Proposition~\ref{prop:keyconstruction} to give a new
proof of Levine's result of~\cite{Levine23} that every strongly negative amphichiral knot bounds a disk in the same rational homology ball, once we show that Construction~\ref{cons:amphi} recovers every strongly negative amphchiral knot. 

\begin{proposition}
    Every strongly negative amphichiral knot bounds a disk in $\Spin(\mathbb{RP}^3)^\circ$. 
\end{proposition}
\begin{proof}
    Let $K$ be a strongly negative amphichiral knot. By the combination of Proposition~\ref{prop:keyconstruction} and Construction~\ref{cons:amphi}, it is enough to show that $K$ has a diagram as in Figure~\ref{fig:rp33}.
    \begin{figure}[h!]
\begin{tikzpicture}
        \node[anchor=south west,inner sep=0] at (0,0)  {\includegraphics[width=0.4\linewidth]{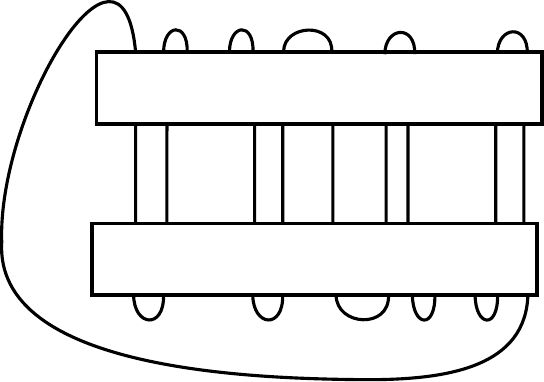}};
         \node at (3.2,1.4){$\text{rot}(\sigma)$};
 \node at (3.2,3.3){$\sigma^{-1}$};
 \end{tikzpicture}
\caption{The local knots obtained from Construction~\ref{cons:amphi},  where $\sigma$ is an arbitrary braid on $4k+1$ strands. } 
    \label{fig:rp33}
\end{figure}
    
    By the resolution of the Smith conjecture, we know that we can realize the involution corresponding to a given SNA knot by rotation about a point on the knot (and the point at infinity, which is on the knot) followed by reflection in the plane of the diagram. We can then equivariantly isotope the knot such that it is in bridge position respecting the involution. In such a position, we know that the braid in the middle must be given by a braid followed by a 180 rotation of the inverse braid.

    A priori, this only gives every SNA knot with odd bridge number, as we must use an even number of $\alpha$- and $\beta$-curves. However, every SNA knot with even bridge number can be put into SNA position with odd bridge number by adding a cancelling pair of 1-dimensional handles, so we can always realize every knot using the above procedure.
\end{proof}

\begin{remark}
    We note that by using a similar cobordism as above in $RP^3 \times I$ instead of $(RP^3)^\circ \times I$, we obtain the fact that every strongly negative amphichiral knot is genus at most 1 in $(RP^3 \times S^1)^\circ$.
\end{remark}

We have also developed a more general method of constructing slice surfaces for Kirby diagrams. A Kirby diagram of a 4-manifold $X$ is a link diagram $D = D_1 \cup D_2$ such that $D_1$ is an unlink,
together with decorations consisting of integers for $D_2$ and dots for $D_1$, which together give instructions for building a 4-manifold (see \cite{GompfStipsicz}). $D_2$ represents the attaching circles for the 2-handles, while $D_1$ represents the 1-handles.

%For example... trefoil (or figure eight? twist knot $5_2$? whatever else we were going to do here.) 

We can build surfaces in such a picture as in $S^3 \times I$, but with an augmented list of moves given by sliding over the various handles.
Using this procedure, we prove Theorem~\ref{thm:mobiusinspinr3}.

\begin{proof}[Proof of Theorem~\ref{thm:mobiusinspinr3}]
We use the handle diagram of Spin($RP^3$) given in Figure \ref{fig:rp3mob}$(iv)$ (ignoring the blue link for now), with one 0-handle, one dotted 1-handle, two 2-handles, and invisible 3- and 4-handles. In this argument, we visualize $\Spin(\mathbb{RP}^3)^\circ$ by removing the 0-handle, and build our surface from the bottom up,  beginning with $K \subset S^3$. 

Our surface begins with a planar cobordism in $S^3 \times I$ from $K$ to the unknot together with some number $k$ of meridians. (One can most easily verify the existence of such a surface  by noting that one can induce a crossing change on a knot by a band move at the cost of adding a meridian to the diagram. We  obtain the link $L$ shown when $k=3$ in Figure~\ref{fig:rp3mob}$(i)$.

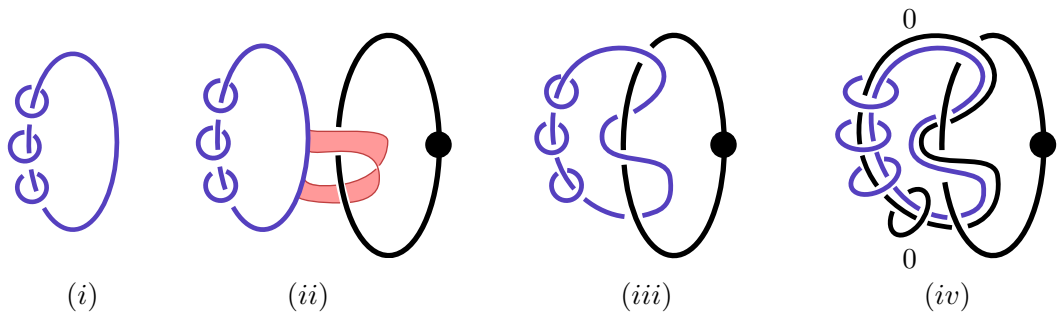
\begin{figure}[h!]
\begin{tikzpicture}
\node[anchor=south west, inner sep=0] at (0,.35){\includegraphics[height=2.4cm]{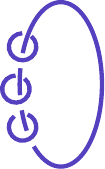}};
\node[anchor=south west, inner sep=0] at (2.5,0){\includegraphics[height=3cm]{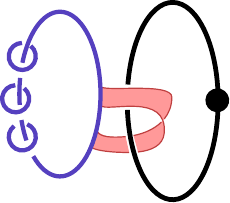}};
\node[anchor=south west, inner sep=0] at (7,0){\includegraphics[height=3cm]{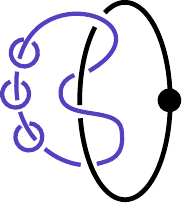}};
\node[anchor=south west, inner sep=0] at (11,0){\includegraphics[height=3cm]{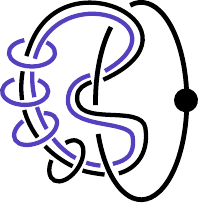}};
\node at (12,3.2){$0$};
\node at (12,0){$0$};
\node at (1, -.5){$(i)$};
\node at (4, -.5){$(ii)$};
\node at (8.5, -.5){$(iii)$};
\node at (12.5, -.5){$(iv)$};
\end{tikzpicture}
\caption{Constructing a surface in $\Spin(\mathbb{RP}^3)^\circ$ with boundary $L$.}
\label{fig:rp3mob}
\end{figure}
We now see the 1-handle of $\Spin(\mathbb{RP}^3)^\circ$ appear, and consider the surface given by attaching an orientation reversing
band to $L$ as shown in Figure~\ref{fig:rp3mob}$(ii)$, resulting in the link $L' \subset S^1 \times S^2$ shown in  Figure~\ref{fig:rp3mob}$(iii)$. We then see the two 2-handles of $\Spin(\mathbb{RP}^3)^\circ$ appear, and observe that the image of $L'$, shown in Figure~\ref{fig:rp3mob}$(iv)$.   bounds $k+1$ disjointly embedded disks, $k$ of which are parallel copies of the core of one of the 2-handles and 1 of which is a core of the other 2-handle. The union of (1) the planar surface, (2) a single orientation-reversing band, and (3) these $k+1$ disks is an embedded M{\"o}bius band with boundary $K$ in $\Spin(\mathbb{RP}^3)^\circ$. 
% \begin{figure}[h!]
% \includegraphics[height=10cm]{image0.pdf}
% \caption{Top left to top right: some number of band moves, each inducing a crossing change at the cost of an added meridian.
% top right to bottom left: a 4-dimensional 1-handle attachment, along with the band to be added after the handle attachment along the handle.
% bottom left to middle left: the band move.
% bottom middle to bottom middle isotopy of the red curve to a writhe 0 diagram so that the framings are clearer.
% bottom middle to bottom right: 2-handle attachments such that all of the red curves are isotopic to appropriately framed attaching circles of the handles.}
%     \label{fig:differenttry}
% \end{figure}
%(maybe add a number of each here, but would need to correspond with labelings in the figure).
%Given a group presentation, we can construct a handlebody that realizes that presentation. By doubling that manifold, we get a closed 4-manifold with the given fundamental group. The process of doubling introduces a 0-framed meridian to each 2-handle, as well as an invisible 3-handle. If we do this process to  $\langle a|a^2\rangle$, we get the spin of $RP^3$.
%Starting with one copy of the relator curve and $n$ of its 0-framed meridians we can unhook the relator from the dotted handle with a single band move. the resulting link is then planar cobordant to any knot of unknotting number at most $n$. (See Figure \ref{fig:rp3mob}). ranging over all $n$, we see that every knot bound a Mobius band by this procedure. 
\end{proof}

% More generally, one can follow a process for finding slice surfaces similar to \cite{Freedman_2010}, namely banding together copies of the 2-handle attaching circles, but with a few key differences. Firstly, the 4-manifolds in question will not be geometrically simply connected, meaning they will have a nonzero number of 1-handles. To obtain a local knot by this procedure, we will need to add bands to unlink the 2-handle attaching circles from the 1-handle curves. To solve this problem, we introduce the other main departure, which is to add genus. By adding more bands than there are 2-handle curves, we can cancel the linking with the 1-handle curves, and obtain a cobordism to a local knot.

% This gives a genus $g$ cobordism in $(S^3 \setminus D_1) \times I$ from appropriately framed longitudes of $D_2$ to some local link $L$. By capping off with cores of the 2-handles, this gives a genus $g$ surface with boundary $L$ in $X^\circ$.

\begin{remark}
      One can use a similar method to find uniform bounds on the nonorientable genus or orientable genus in the double of 1-2 handlebodies. We briefly describe the process, leaving the details to the interested reader. 
       Let $\langle g_1, \dots g_n|r_1, \dots r_m\rangle$ be a group presentation, and let $D$ be the corresponding Kirby diagram of the double of the 1-2 handlebody with that group presentation. Let $R$ be a subset of the relator set such that every generator appears an even number of times in the words of the relators. Using a 2-handle core for each element of $R$, we obtain a link that links each 1-handle dotted curve an even number of times. By adding (potentially orientation reversing) band moves, we can cancel this linking in pairs as in Theorem \ref{thm:mobiusinspinr3}, to obtain a local link. 
       %Using the 0-framed meridians of each relator curve, we see that this will give a uniform bound on the genus. We note this process is equivalent to finding an embedded surface in your manifold with a 0-framed sphere dual.
For example, the spin of $L(4,1)$ is the double of the 1-2 handlebody corresponding to the presentation $\langle a|a^4\rangle$. 
This will require 2 orientation reversing bands to unlink the 2-handle from the 1-handle, meaning all knots will have nonorientable genus at most 2 in this manifold. However, the double of the presentation $\langle a,b|a^2, ab^2\rangle$ also gives spin of $L(4,1)$, and doing the same process to only the 2nd relator curve gives an upper bound on the nonorientable genus of 1.  
%Doubling the presentation for  $\langle a,b|aba^{-1}b^{-1}\rangle$, we see that we also need 2 bands to unlink the 2-handle, but in this case both bands can be made orientation preserving, so we find that all knots are orientable genus 1 in this manifold. This coincides with our intuition, as this manifold is $T^2 \times S^2$, which has an embedded torus with sphere dual.
\end{remark}

%(must have a collection of curves that links all generators an even number of times for non ori, must be algebraically 0 for ori)

%\begin{example}
%--SNA bound in Spin(RP3) \\
%--SNA are genus in 1 in RP3 x S1
%\end{example}

%\subsection{Other stuff}

%\begin{example}
 %   A knot bounding a small genus surface in a ZHB4. Unclear if I can come up with a reasonable one of these because of the monodromy.
%\end{example}

%For more general manifolds (i.e. homology spheres, homotopy spheres), we can still get interesting families of knots and links using cancelling pairs of relator curves and banding them, but in this case the meridians induce band passes instead of crossing changes (i.e.passing a strand through 2 strands at once in certain spots)

\subsection{$T^4$ genus and $S^4$ genus differ}~\label{sec:t4}

We now detour slightly to consider the difference between the minimal genus of a surface bounded in $(T^4)^{\circ}$ versus $B^4$ for a given knot. 
 By taking covers, one can immediately see that a knot bounds a disk in $(T^4)^\circ$ if and only if it is smoothly slice. In fact, the following more general proposition holds. 

\begin{proposition}~\label{prop:genuscollapset4}
    Let $K \subset S^3$ be a knot, and let $F \subset (T^4)^\circ$ be a genus $g$ surface  with boundary $K$. Suppose that the inclusion induced map $i_*: \pi_1(F) \to \pi_1((T^4)^\circ) = \mathbb{Z}^4$ is not full rank, i.e. there is a primitive element $x \in \pi_1((T^4)^\circ)$ such that no  power of $x$ is in $im(i_*)$. 
    Then $K$ bounds a genus $g$ surface in $B^4$.
\end{proposition}
\begin{proof}
    Take the cover of $(T^4)^\circ$ corresponding to the subgroup generated by $x$. Because $x$ is a primitive element, this cover is $T^3 \times\mathbb{R}$ punctured infinitely many times. Because $x$ is not in $im(i_*)$, $F$ lifts to infinitely many disjoint genus-$g$ surfaces. To complete the proof, we note that $T^3 \times\mathbb{R}$ embeds in $S^4$, for example as the tubular neighborhood of the boundary of a torus.
\end{proof}

In particular, naive attempts to apply an analogue of Proposition~\ref{prop:keyconstruction} to $T^4=T^3 \times S^1$ are doomed to fail: one can check that the resulting surfaces will always satisfy the hypotheses of Proposition~\ref{prop:genuscollapset4}. 
We also obtain the following corollary. 

\begin{corollary}
For $g \in \{0,1\}$, $K$ bounds an orientable  surface of genus $g$ in $(T^4)^{\circ}$ only if it bounds such a surface in $B^4$.

For $h \in \{1,2,3\}$, $K$ bounds a nonorientable surface with first Betti number $h$ in $(T^4)^{\circ}$ only if it bounds such a surface in $B^4$.
\end{corollary}

We are nonetheless able to prove the following, which immediately implies Theorem~\ref{thm:g4t4}.
\begin{proposition}
For $m$ odd, let $L_m$ be the $2m$-component link obtained by replacing both components of the Hopf link by $m$ coherently oriented 0-framed parallel copies, and let $K$ be any knot obtained by attaching $2m-1$ bands to $L$. 
Then $K$ bounds a surface in $(T^4)^\circ$ of genus $2m$, even though $g_4(K) \geq m^2/2-m+1$.
\end{proposition}

\begin{proof}
We first argue that $L_m$ bounds a disjoint union of $2m$ tori in $(T^4)^\circ$. Considering $T^4= T^2 \times T^2$, we see that $T^2 \times \{*\}$ and $\{*\} \times T^2$ form 2 0-framed tori that meet in a single transverse intersection. Take $m$ parallel copies of each,  and remove a small ball around the intersection point to see $L_m$ on the boundary. 
Therefore, any knot $K$ obtained by attaching $2m-1$ bands to $L_m$ bounds a surface of genus $2m$ in $(T^4)^\circ$. 

Now, we need only compute that $\sigma(L_m)=-m^2$. This will imply that as desired
\[ g_4(K) \geq \frac{1}{2} |\sigma(K)| \geq \frac{1}{2} \left|-m^2+2(m-1)\right| \]
since the signature function behaves predictably under banding~\cite{Tristram}. 

Proposition 2.5 of Cimasoni-Florens~\cite{CimasoniFlorens} implies that
 \[ \sigma(L_m)= \sigma_{\text{col}}(L_m; -1,\dots, -1, -1,\dots,-1)-\text{lk}(L_m),\]
where $\sigma_{\text{col}}(L_m; -1,\dots, -1, -1,\dots,-1)$ denotes the signature of $L_m$ considered as a $2m$-colored link, and $\text{lk}(L_m)$ denotes the total linking number of $L_m$, i.e. the sum of the linking numbers of all 2-component sublinks of $L_m$. Since $\text{lk}(L_m)=m^2$, it suffices to show that $\sigma_{\text{col}}(L_m; -1, -1, \dots, -1)=0$.
By Theorem 2.10 of~\cite{DFL17}, we have that 

\begin{align*}
\sigma_{\text{col}}(L_m; -1, \dots, -1)&=
\left(\text{ind}\left(\sum_{i=1}^m \text{Log}(-1) \right)-\sum_{i=1}^m
\text{ind}(\text{Log}(-1))\right)^2\\
&= \left(\text{ind}\left(\frac{m}{2}\right)- \sum_{i=1}^m \text{ind}\left(\frac{1}{2}\right)\right)^2\\
&= \left(\left\lfloor \frac{m}{2} \right\rfloor -\left\lfloor \frac{-m}{2} \right\rfloor -m\left(\left\lfloor \frac{1}{2} \right\rfloor -\left\lfloor \frac{-1}{2} \right\rfloor \right)\right)^2\\
&=0.
\end{align*}
(Here $\text{Log} \colon S^1 \to [0,1)$ is characterized by the property that $\text{Log}(e^{2 \pi i s})=s$ for all $s \in [0,1)$ and $\text{ind} \colon \mathbb{R} \to \Z$ is defined by $\text{ind}(x)= \lfloor x \rfloor - \lfloor -x \rfloor$.)
%Note that for any integer $n$, $\text{ind(n/2)}=n$.
\end{proof}

\section{Knots bounding small genus surfaces in $\Spin(L(3,1))^\circ$}\label{sec:CGstuff}

We now begin to define the knots featuring in Theorems~\ref{thm:mainthm} and~\ref{thm:biggaps}. 

\begin{definition}~\label{defn:ll7strand}
For a pure 7-strand string link $\mathcal{L}$, let $L_{\mathcal{L}}$ be as in Figure~\ref{fig:LL}, where $\mathcal{L}^{-1}$ denotes the inverse of $\mathcal{L}$ in the string link concordance group~\cite{StringLinks}. 

\begin{figure}[h!]
\begin{tikzpicture}
%\draw[step=1cm,color=gray] (0,0) grid (3,5);
\node[anchor=south west, inner sep=0] at (0,1){\includegraphics[height=5cm]{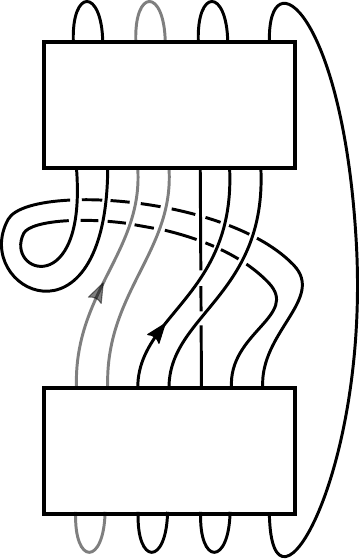}};
\node at (1.6,2){$\mathcal{L}^{-1}$};
\node at (1.5,5){$\mathcal{L}$};
\end{tikzpicture}
\caption{Given a pure strink link  $\mathcal{L}$, the link $L_{\mathcal{L}}$ is a 2-component link. }
\label{fig:LL}
\end{figure}
\end{definition}

\begin{proposition}~\label{prop:cobbuildingl31}
    For any pure 7-strand string link $\mathcal{L}$, the 2-component link $L_{\mathcal{L}}$ bounds a planar surface in $\Spin(L(3,1))^\circ$.
    \end{proposition}
\begin{proof}
We will construct a planar cobordism $\Sigma$ between $T_J \# L_{\mathcal{L}}$ and $T_J$ in $L(3,1)^{\circ} \times I$ for some knot $J \subset L(3,1)$, in order to apply Proposition~\ref{prop:keyconstruction}. 
As before, let $H$ be the usual Heegaard torus for $L(3,1)$, and fix a Heegaard diagram for $L(3,1)$ where the $\alpha$-curve is the $(1,3)$ curve and the $\beta$-curve is the $(1,0)$ curve. Let $J$ be the image of the $(0,1)$ curve in $L(3,1)$. 
The only part of our cobordism  that is not contained in $\nu(H)^{\circ} \times I$ will be disks lying in the handlebodies of the Heegaard splitting with boundaries that are copies of the $\alpha$- and $\beta$- curves of our fixed Heegaard diagram, attached at $t=1/4$. 

For $I \subset [0,1]$, we let $\Sigma_I:=\Sigma \cap (L(3,1)^\circ\times I)$, and define

\begin{enumerate}
\item $\Sigma_{[0,\frac{1}{4}]}$ as $T_J \times [0,1/4]$ together with two oppositely oriented $\alpha$-disks and two oppositely oriented $\beta$-disks attached at $t=1/4$, resulting in $\Sigma_{\frac{1}{4}}$  as shown on the left of Figure~\ref{fig:l31firstbit}. 

\item $\Sigma_{[\frac{1}{4},\frac{1}{2}]}$  as a string link concordance from the  trivial 7-strand string link to $\mathcal{L}^{-1} * \mathcal{L}$, resulting in $\Sigma_{\frac{1}{2}}$ as shown in the right of Figure~\ref{fig:l31firstbit}, ignoring the orange bands. 

\begin{figure}[h]
    \begin{tikzpicture}
        \node[anchor=south west,inner sep=0] at (0,0)
{\includegraphics[height=4.25cm]{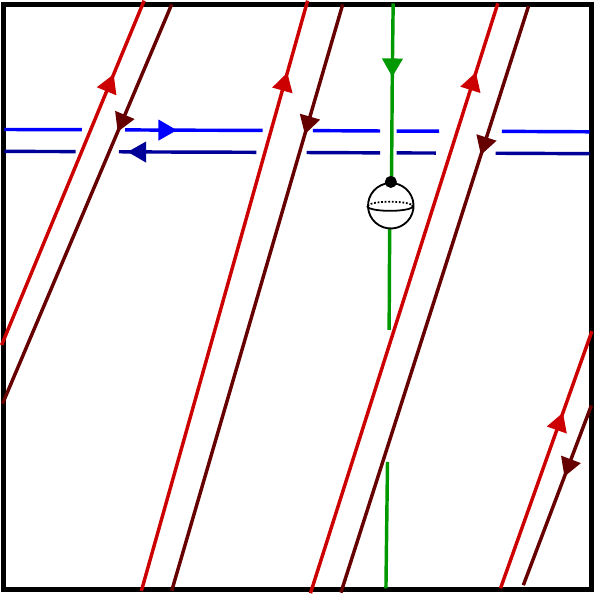}};
        \end{tikzpicture}
        \qquad  
    \begin{tikzpicture}
        \node[anchor=south west,inner sep=0] at (0,0)
{\includegraphics[height=4.25cm]{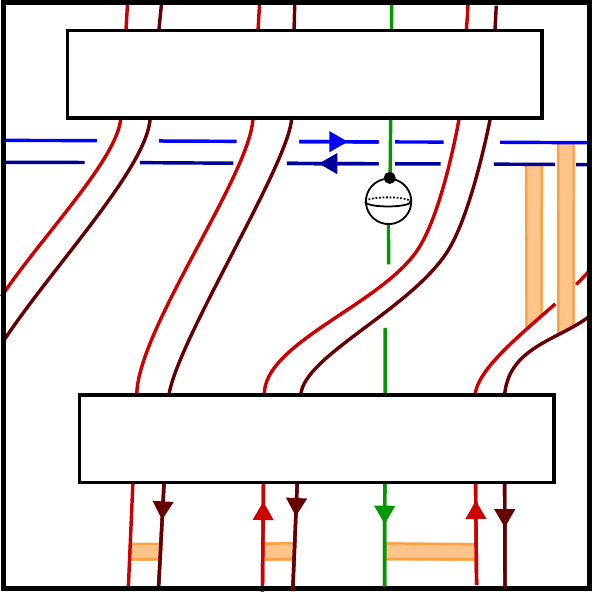}};
\node at (2.2,3.75){$\mathcal{L}$};
\node at (2.4,1.15){$\mathcal{L}^{-1}$};
        \end{tikzpicture}
        \qquad
    \caption{The part of $\Sigma$ in $\nu(H)^\circ \times \{t\}$ for $t=\frac{1}{4}$ (left) and $\frac{1}{2}$ (right).}
    \label{fig:l31firstbit}
\end{figure}

\item $\Sigma_{[\frac{1}{2},\frac{3}{4}]}$  as obtained from $\Sigma_{\frac{1}{2}}$ by attaching 5 oriented orange bands, resulting in $\Sigma_{\frac{3}{4}}$ as shown on the right of Figure~\ref{fig:l31firstbit}.

\item $\Sigma_{[\frac{3}{4},1]}$ is obtained as described in Example~\ref{exl:fig8!} by point-pushing along the curve indicated on the left of Figure~\ref{fig:pointpush}, resulting in $\Sigma_1=T_J \#L_{\mathcal{L}}$ as shown on the right of Figure~\ref{fig:pointpush}.
\end{enumerate}

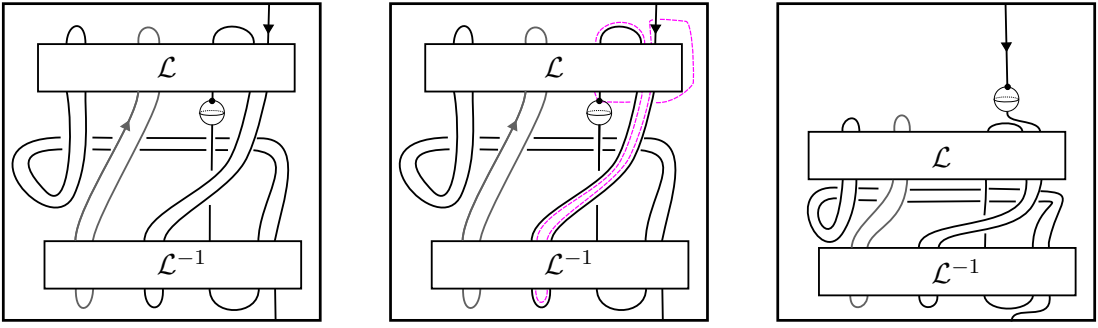
\begin{figure}[h]
\begin{tikzpicture}
        \node[anchor=south west,inner sep=0] at (0,0)
{\includegraphics[height=4.25cm]{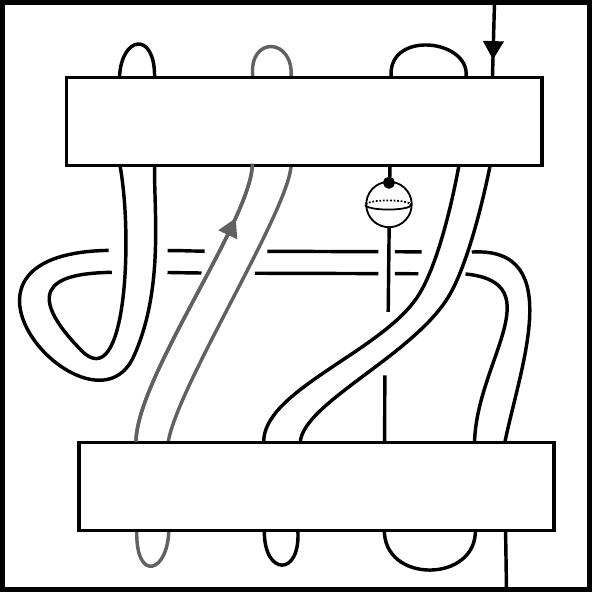}};
\node at (2.2,3.4){$\mathcal{L}$};
\node at (2.4,.8){$\mathcal{L}^{-1}$};
\end{tikzpicture}
\qquad
\begin{tikzpicture}
        \node[anchor=south west,inner sep=0] at (0,0)
{\includegraphics[height=4.25cm]{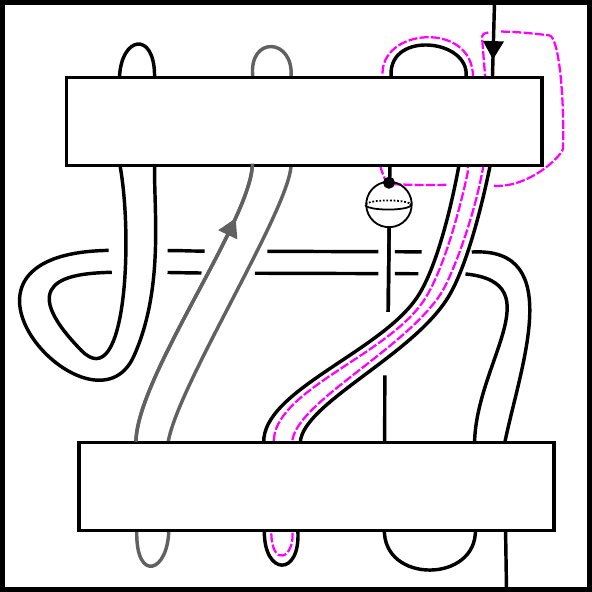}};
\node at (2.2,3.4){$\mathcal{L}$};
\node at (2.4,.8){$\mathcal{L}^{-1}$};
\end{tikzpicture}
\qquad 
\begin{tikzpicture}
        \node[anchor=south west,inner sep=0] at (0,0)
{\includegraphics[height=4.25cm]{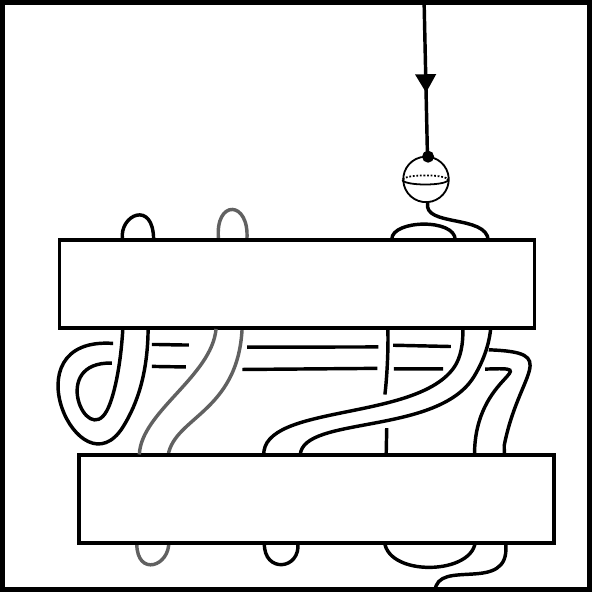}};
\node at (2.2,2.2){$\mathcal{L}$};
\node at (2.4,.68){$\mathcal{L}^{-1}$};
\end{tikzpicture}
\caption{The part of $\Sigma$ in $\nu(H)^\circ \times \{t\}$ for $t=\frac{1}{2}$ (left), $\frac{3}{4}$ (center), and $1$ (right).}\label{fig:pointpush}
\end{figure}
Applying Proposition~\ref{prop:keyconstruction} now finishes the proof. 
\end{proof}

\begin{remark}
We will work with knots $K$ obtained from $L_{\mathcal{L}}$ by a single orientation preserving band move. Such knots bound a genus one surface in $\Spin(L(3,1))^{\circ}$, but do not obviously bound a genus one surface in $B^4$. (One can check that they do bound genus two surfaces in $B^4$.)
By construction, such a knot $K$  has the property that $T_J \# K$ and $T_J$ cobound a genus one surface in $L(3,1)^{\circ} \times I$.  By gluing on an unknotted $(B^3,I) \times I$ pair, we see that $J \#K$ and $J$ cobound a genus one surface $F$ in $L(3,1) \times I$. When we take the universal cover, $F$ lifts to $\widetilde{F}$, a genus 3 surface cobounded by $\widetilde{J\#K}= K \#K \# K$ and $\widetilde{J}=U$.  We therefore have that $g_4(K\#K\#K)\leq 3$ and hence, as noted in the introduction, that many 4-genus bounds are incapable of showing that $g_4(K)>1$. 

For the specific examples of Theorem~\ref{thm:mainthm}, we know that $g_4(K)=2$ and $g_4(K\# K \# K) \leq 3$. It is interesting to ask what the precise values of $g_4(nK)$ are for $n \geq 2$, or even what the value of the stable 4-genus $g_{st}(K):=\lim_{n\to \infty} \frac{g_4(nK)}{n}$ of~\cite{LivingstonStable}  might be.  (It is certainly no more than 1, and it is not clear that it is nonzero.) Unhelpfully, for our examples the Levine-Tristram signature function of $K$ is identically zero:  one can compute $\Delta_K(t)$ and observe that it has no roots on the unit circle $S^1 \subset \mathbb{C}$. 
\end{remark}

\begin{remark}[Interesting 3-manifolds from interesting annuli]

Given a pure 7-strand string link $\mathcal{L}$, let $\mathcal{A}_{\mathcal{L}}$ be the annulus in $\Spin(L(3,1))^\circ$ bounded by $L_{\mathcal{L}}$ that is obtained by combining the construction above with Proposition~\ref{prop:keyconstruction}. By appropriately choosing $\mathcal{L}$, we can ensure that the linking number of $L_\mathcal{L}$ is zero. Doing $1/n$ surgery on $\Spin(L(3,1))^\circ$ along $\mathcal{A}_{\mathcal{L}}$ will then give integer homology 3-spheres $Y_{\mathcal{L}}$ that by construction bound rational homology 4-balls, but which have no particular reason to bound integer homology 4-balls. 
 Notably, the rational homology balls obtained from this process are $\mathbb{Z}/2\mathbb{Z}$ homology balls, so many of the invariants one might otherwise use in such a situation (the Rokhlin invariant, involutive Heegaard Floer homology, etc), will vanish. So obstructing $Y_{\mathcal{L}}$ from bounding an integer homology 4-ball may be beyond the scope of current obstructive techniques.

\end{remark}

\begin{example}
For an ordered tuple of integers $\mathbf{n}=(n_0,n_1,n_2,n_3,n_4)$ and an ordered tuple of knots $\mathbf{J}=(J_0, J_1, J_2,J_3,J_4)$,  let $\mathcal{L}= \mathcal{L}(\mathbf{n},\mathbf{J})$ denote the corresponding 7-strand string link illustrated in Figure~\ref{fig:stringlink}. (Note that $n_0$ and $J_0$ do not contribute to $\mathcal{L}_{\mathbf{n}, \mathbf{J}}$, but will play a role later on.)
%\footnote{AM: $n_0$, $J_0$ do not contribute to $\mathcal{L}_{\mathbf{n}, \mathbf{J}}$}.
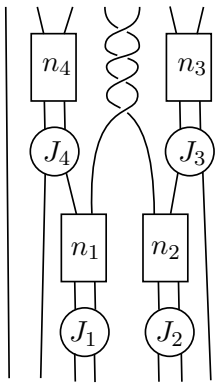
\begin{figure}[h!]
\begin{tikzpicture}
%\draw[step=1cm,color=gray] (0,0) grid (3,5);
%Uncomment this to get some helpful grid lines
\node[anchor=south west, inner sep=0] at (0,0){\includegraphics[height=5cm]{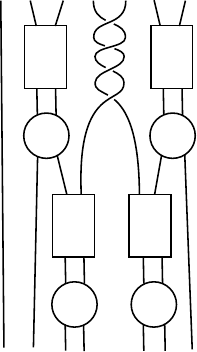}};
\node at (1.05,.65){$J_1$}; \node at (1.08,1.75){$n_1$};
\node at (2.2,.65){$J_2$};\node at (2.15,1.76){$n_2$};
\node at (.7,3.05){$J_4$}; \node at (.7,4.2){$n_4$};
\node at (2.5,3.05){$J_3$}; \node at (2.5,4.2){$n_3$};
\end{tikzpicture}
%a=n_4, b=n_3, c=n_1, d=n_2, m=n_0
\caption{A 7-strand pure string link $\mathcal{L}(\mathbf{n},\mathbf{J})$, where the rectangular boxes denote full twists and the circular regions depict the result of tying the two parallel strands into the given knot.}
\label{fig:stringlink}
\end{figure}
\end{example}

\begin{figure}[h!]
\begin{tikzpicture}
%\draw[step=1cm,color=gray] (0,0) grid (14,11);
\node[anchor=south west, inner sep=0] at (0,0){\includegraphics[height=11cm]{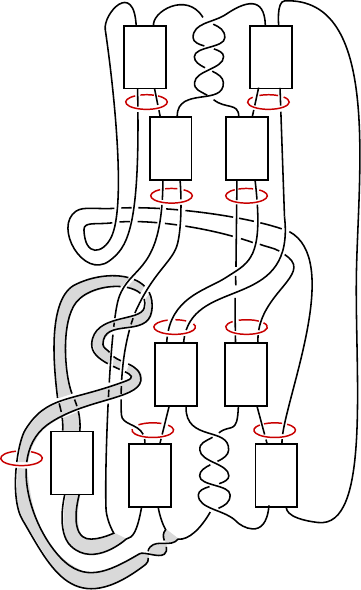}};
%a=n_4, b=n_3, c=n_1, d=n_2, m=n_0
\node at (2.7,9.85){\small $n_4$}; \node at (5.1,9.85){\small $n_3$};
\node at (3.3, 9.25){\small $\gamma_4$};
\node at (5.6, 9.1){\small $\gamma_3$};
\node at (3.15,8.2){\small $n_1$}; \node at (4.6,8.2){\small $n_2$};
\node at (3.8,7.5){\small $\gamma_1$}; \node at (5.05,7.1){\small $\gamma_2$};
\node at (3.85,4.8){\small $\gamma_1'$};\node at (5.2,4.8){\small $\gamma_2'$};
\node at (3.3,4){\small $-n_1$}; \node at (4.57,4){\small $-n_2$};
\node at (.7,2.1){\small $\gamma_0$}; \node at (1.4,2.3){\small $n_0$}; \node at (3.4,2.8){\small $\gamma_4'$};
\node at (2.8,2.1){\small $-n_4$}; \node at (5.13,2.13){\small $-n_3$}; \node at (5.7,2.8){\small $\gamma_3'$};
\node[anchor=south west, inner sep=0] at (7,0){\includegraphics[height=11cm]{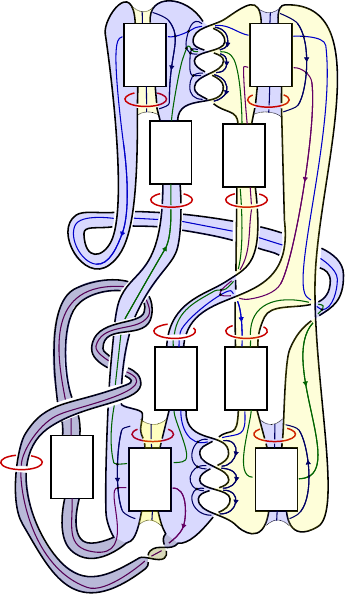}};
\node at (9.7,10){\small $n_4$}; \node at (12,10){\small $n_3$};
\node at (10.18,8.1){\small $n_1$}; \node at (11.52,8.1){\small $n_2$};
\node at (10.3,4){\small $-n_1$}; \node at (11.55,4){\small $-n_2$};
\node at (8.4,2.3){\small $n_0$}; 
\node at (9.8,2.1){\small $-n_4$}; \node at (12.1,2.13){\small $-n_3$}; 
\node at (10.3,10.65){\tiny$_0$}; \node at (11.35,10.2){\tiny$_1$}; 
\node at (12.8,9.9){\tiny$_2$};\node at (9.4,8.5){\tiny$_3$};\node at (10.05,6.2){\tiny$_4$};
 \node at (12.65,8.5){\tiny$_5$}; 
\node at (9.3,3.2){\tiny$_6$};
\node at (11.45, 2.5){\tiny$_7$};
\node at (12.8,1.7){\tiny$_8$};
\node at (10.35,1.5){\tiny$_9$};
\node at (11.6, 5.2){\tiny$_{10}$};
\node at (12.6,3.5){\tiny$_{11}$};
\node at (11.28,9.65){\tiny$_{12}$};\node at (11.3,9.2){\tiny$_{13}$};
\node at (11.5, 2){\tiny$_{14}$};
\node at (11.5, 1.6){\tiny$_{15}$};
\end{tikzpicture}

\caption{The base knot $K(\mathbf{n})$ (left) and a genus eight Seifert surface $F$ for $K(\mathbf{n})$ with a collection of curves $x_0, \dots, x_{15}$ representing a basis for $H_1(F)$(right). Note that on the right side of the figure,  $x_i$ is labeled by $i$ and we omit the labels of the $\gamma$-curves.}
\label{fig:baseknot}
\end{figure}
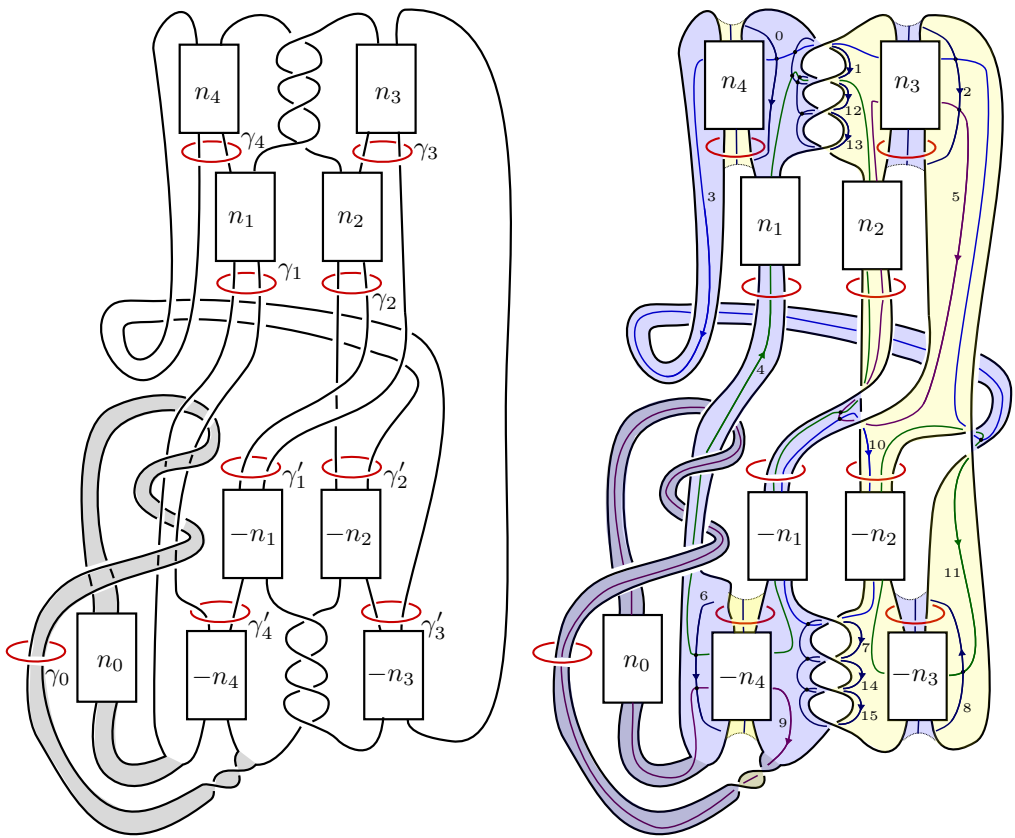

Let $K(\mathbf{n})$ be as on the left of Figure~\ref{fig:baseknot}. 
For a tuple $\mathbf{J}$ as above,  let 
$K(\mathbf{n})_{\mathbf{J}}$ denote the result of infecting $K(\mathbf{n})$ by $J_0$ along $\gamma_0$,  by $J_i$ along $\gamma_i$ for $i=1,2,3,4$ and by $-J_i$ along $\gamma_i'$ for $i=1,2,3,4$.  
Observe that $K(\mathbf{n})_{\mathbf{J}}$ is obtained from $L_{\mathcal{L}(\mathbf{n},\mathbf{J})}$ as in Definition~\ref{defn:ll7strand} by adding a band. 
Therefore, $K(\mathbf{n})_{\mathbf{J}}$
 bounds a genus one surface in $\text{Spin}(L(3,1))^\circ$, and so Theorem~\ref{thm:mainthm} is an immediate consequence of the next result. 

\begin{proposition}\label{prop:g4is2}
Let $K=K(-2,-2,2,1,-1)$.  
%a=-1,b=1,c=-2,d=2,m=-2
%a=n_4, b=n_3, c=n_1, d=n_2, m=n_0
There exists $\mathbf{J}$ such that $g_4(K_{\mathbf{J}}) \geq 2$.  
\end{proposition}
%Note that the Alexander polynomial of K is (-2 + x)^2 * (-1 + 2*x)^2 * (2 - 6*x + 8*x^2 - 14*x^3 + 21*x^4 - 14*x^5 + 8*x^6 - 6*x^7 + 2*x^8), which one can check has no roots on the unit circle. So K has trivial Levine-Tristram signature function. 

In the proof we will give explicit conditions on $\mathbf{J}$ ensuring that $g_4(K_{\mathbf{J}}) \geq 2$ in terms of the Tristram-Levine signature functions of $J_0,J_1,J_2,J_3,J_4$.

Without going into the details of the definition, we remind the reader that Casson and Gordon associate a rational number $\sigma(K,\chi)$ to a knot $K$ and a homomorphism $\chi \colon H_1(\Sigma_K) \to \Z/q \Z$ by considering signatures of certain twisted intersection forms, see~\cite{CassonGordon2}.
To prove Proposition~\ref{prop:g4is2}, we will imitate the proof of~\cite[Proposition 2.7]{MillerAmphi} as follows. 
Roughly speaking, Gilmer's Theorem~\ref{thm:Gilmer} tells us that if $g_4(K_{\mathbf{J}})=1$, then at least nine Casson-Gordon signatures are ``small''. Simultaneously, Litherland's Theorem~\ref{thm:Litherland} tells us that the Casson-Gordon signatures of $K_{\mathbf{J}}$ are controlled by the Tristram-Levine signatures of the knots in $\mathbf{J}$, together with the value of $\chi$ on the infection curves $\gamma_i, \gamma_i'$. In particular, by our choice of $\mathbf{J}$ we can ensure that the Casson-Gordon signatures of $K_{\mathbf{J}}$ are ``large'' unless $J_0$ makes no contribution and $J_i$ and $-J_i$ make canceling contributions for $1 \leq i \leq 4$. 
The technical heart of the proof is Lemma~\ref{lemma:charactertracking}, which shows that there are only three $\chi$ for which this occurs, and hence only three characters which could have ``small'' Casson-Gordon signatures, thereby contradicting Gilmer.  

We begin by stating a special case of Litherland's formula for the Casson-Gordon signatures of satellite knots, which we will apply to $K_{\mathbf{J}}$.
\begin{theorem}[\cite{Litherland}]~\label{thm:Litherland}
Let $P$ be a winding number 0 pattern described by a curve $\eta$ in the complement of $P(U)$ and let $J$ be a knot. 
There is a canonical isomorphism 
$ \alpha \colon H_1(\Sigma_{P(J)}) \to H_1(\Sigma_{P(U)})$ with the property that
for $q$  a prime power and $\chi \colon H_1(\Sigma_{P(U)}) \to \Z/q\Z$,
\[ \sigma(P(J), \chi \circ \alpha)= 
\sigma(P(U), \chi) + 2 \sigma_J(\omega_q^{\chi(\widetilde{\eta})}),
\]
where $\widetilde{\eta}$ denotes a lift of $\eta$ to $\Sigma_{P(U)}$ and  $\omega_q= e^{2 \pi i/ q}$. 
\end{theorem}

 \begin{remark}
     Litherland's formula involves a sum of $ \sigma_J(\omega_q^{\chi(\widetilde{\eta}_j)})$ across the 
 $n$ distinct lifts of $\eta$ to the $n$-fold cyclic branched cover of $P(K)$.  However, when $n=2$ the homology classes of the lifts $\widetilde{\eta}_1$ and $\widetilde{\eta}_2$ satisfy $[\widetilde{\eta}_1]= - [\widetilde{\eta}_2]$,  and so by the symmetry of the Tristram-Levine signature function we have
 \[ \sigma_J(\omega_q^{\chi(\widetilde{\eta}_1)})+ \sigma_J(\omega_q^{\chi(\widetilde{\eta}_2)})=
 \sigma_J(\omega_q^{\chi(\widetilde{\eta}_1)})+ \sigma_J(\omega_q^{-\chi(\widetilde{\eta}_1)})
 =2 \sigma_J(\omega_q^{\chi(\widetilde{\eta}_1)}),\]
thereby obtaining the statement of Theorem~\ref{thm:Litherland}. 
 \end{remark}

We will henceforth  abuse notation by suppressing $\alpha$ and failing to  distinguish between $\chi \colon H_1(\Sigma_{K}) \to \Z_q$ and $\chi \circ \alpha \colon H_1(\Sigma_{K_{\mathbf{J}}}) \to \Z_q$. 
Applying this to $K_{\mathbf{J}}$,  which is obtained from $K$ via nine winding number zero satellite operations,  we have that for any character $\chi \colon H_1(\Sigma_K) \to \Z/ q \Z$, 
\begin{align}\label{eqn:CGsigformula}
 \sigma(K_{\mathbf{J}}, \chi)= 
\sigma(K,  \chi) + 
2 \sigma_{J_0}(\omega_q^{\chi(z_0)})
+ 2 \sum_{i=1}^4 \left(\sigma_{J_i}(\omega_q^{\chi(z_i)})- \sigma_{J_i}(\omega_q^{\chi(z_i')})\right),
\end{align}
where $z_i:=[\widetilde{\gamma_i}]$ for $i=0, \dots, 4$ and  $z_i':= [\widetilde{\gamma_i}']$ for $i=1, \dots, 4$. 
%we use that $\sigma_{-J_i}(\omega)= -\sigma_{J_i}(\omega)$ for all $\omega \in S^1$. 
%Note that we can choose $J_0, J_1, J_2,J_3,J_4$ to have Tristram Levine signature functions so that Equation~\ref{eqn:CGsigformula} is ``large" unless $\chi(z_0)=0$ and $\chi(z_i')= \pm \chi(z_i)$ for $i=1, \dots, 4$. 
%\footnote{All we need   is $\sigma_{J_i}(\omega_9)$, $\sigma_{J_i}(\omega_9^2)$, $\sigma_{J_i}(\omega_9^3)$, $\sigma_{J_i}(\omega_9^4)$ to be increasingly large as we increase $i$ and also very different from each other.  Taking each $J_i$ to be a sufficiently large connected sum of $T(2,9)$s should do the trick,  but also there are results telling us we can realize basically any TL signature function we want.}
Note that if $\chi(z_0)=0$ then $J_0$ makes no contribution to $ \sigma(K_{\mathbf{J}}, \chi)$. Similarly, for $1 \leq i \leq 4$ if $\chi(z_i)= \pm \chi(z_i')$ then $J_i$ and $-J_i$ make a net contribution of zero to  $ \sigma_1(K_{\mathbf{J}}, \chi)$.  We now analyze when this occurs.

\begin{lemma}\label{lemma:charactertracking}
Let $K=K(-2,-2,2,1,-1)$, and let $z_0, z_1, \dots, z_4, z_1', \dots, z_4'$ denote the homology classes of lifts of $\gamma_0, \gamma_1,\dots, \gamma_4, \gamma_1', \dots, \gamma_4'$ to $\Sigma_K$.  
Then  
\begin{enumerate}
\item $H_1(\Sigma_K) \cong (\Z/9 \Z)^4$
\item Up to re-scaling,  there is exactly one non-trivial map $\chi \colon H_1(\Sigma_K) \to \Z/ 3 \Z$ with the property that $\chi(z_0)=0$ and $\chi(z_i)= \pm \chi(z_i')$ for $i=1,2,3,4$.
\item There is no surjective map $\chi \colon H_1(\Sigma_K) \to \Z/9 \Z$ with the property that $\chi(z_0)=0$ and $\chi(z_i)= \pm \chi(z_i')$ for $i=1,2,3,4$.
\end{enumerate}
\end{lemma}

In the proof of Lemma~\ref{lemma:charactertracking}, we will need the following standard result, for which we follow the notation of \cite{FriedlPowell_Blanchfield}; see also~\cite{AMMMPSBranchedCover} for a similar application. 
\begin{proposition}
Let $F$ be a Seifert surface for a knot $K$ with a collection of simple closed curves $x_0,x_1, \dots, x_{2g-1}$ that form a basis for $H_1(F; \Z)$. 
Let $\hat{x}_0,  \hat{x}_1, \dots,  \hat{x}_{2g-1}$ be a collection of simple closed curves in $S^3 \setminus \nu(F)$ representing a basis for $H_1(S^3 \setminus \nu(F); \Z)$ satisfying $\text{lk}(x_i, \hat{x}_j)= \delta_{i,j}$. Let $y_0,y_1,\dots, y_{2g-1}$ denote (the homology classes of) lifts of $\hat{x}_0,  \hat{x}_1, \dots,  \hat{x}_{2g-1}$ to a single copy of $S^3 \setminus \nu(F)$ within
\[ \Sigma_K= X_2(K) \cup S^1 \times D^2=  (S^3 \setminus \nu(F)) \cup (S^3 \setminus \nu(F)) \cup (S^1 \times D^2).\]

Then the first homology of $\Sigma_K$ is generated by $y_0, y_1, \dots, y_{2g-1}$,  with relations given by the columns of $A+A^T$, where $A=[\text{lk}(x_i, x_j^+)]$ is the Seifert matrix. 
\end{proposition}
%We will apply the previous proposition in a way that mimics arguments from~\cite{AMMMPSBranchedCover}. 
\begin{proof}[Proof of Lemma~\ref{lemma:charactertracking}]
Figure~\ref{fig:baseknot} shows a genus eight Seifert surface $F$ for $K$ and a basis $x_0, \dots, x_{15}$ for $H_1(F)$ with respect to which we have the following Seifert matrix:
%a=n_4, b=n_3, c=n_1, d=n_2, m=n_0
\begin{align*} 
A= \left[\begin{array}{ccc|ccc|ccc|ccc|cc|cc}
-2& \cdot & \cdot & 1 & \cdot & \cdot & \cdot & \cdot & \cdot & \cdot & \cdot & \cdot &\cdot  & \cdot & \cdot  &\cdot \\  
\cdot & -1 & \cdot & 0 & 1 & \cdot & \cdot  & \cdot  & \cdot  & \cdot  & \cdot  & \cdot &1  & \cdot & \cdot  &\cdot  \\
 \cdot & \cdot & 1 & 0  & \cdot & 0  & \cdot & \cdot & \cdot & \cdot & \cdot  & \cdot &\cdot  & \cdot & \cdot  &\cdot  \\ \hline
 0 & 1 & 1 & 0   & \cdot & \cdot & \cdot & \cdot & \cdot& \cdot & \cdot & -1 &\cdot  & \cdot & \cdot  &\cdot\\
 \cdot & 0  & \cdot  & \cdot & 2 & -2 & 0  & \cdot  & \cdot  & \cdot & -2 & \cdot &-1  & \cdot & \cdot  &\cdot \\
 \cdot  & \cdot & -1  & \cdot & -3 & 2  & \cdot & \cdot & \cdot & \cdot & 1 & \cdot &\cdot  & \cdot & \cdot  &\cdot\\ \hline
 \cdot  & \cdot  & \cdot  & \cdot & 1  & \cdot & 1  & \cdot  & \cdot  & -1 & \cdot  & \cdot &\cdot  & \cdot & \cdot  &\cdot \\
 \cdot  & \cdot & \cdot & \cdot & \cdot & \cdot & \cdot & 1  & \cdot  & \cdot & -1 & \cdot &\cdot  & \cdot & 0  &\cdot\\
 \cdot  & \cdot & \cdot & \cdot & \cdot & \cdot & \cdot & \cdot & -1 & \cdot & \cdot & 0 &\cdot  & \cdot & \cdot  &\cdot \\ \hline
 \cdot  &\cdot  & \cdot &\cdot  &\cdot  & \cdot & 0 & \cdot & \cdot &-2& \cdot & \cdot &\cdot  & \cdot & \cdot  &\cdot \\
 \cdot  &\cdot  & \cdot  & \cdot & -2  &0  &\cdot  &0  & \cdot & \cdot & 0 & 2 &\cdot  & \cdot & \cdot  &\cdot\\
\cdot  &\cdot  & \cdot & 0& \cdot & \cdot & \cdot & \cdot & -1 & \cdot & 2 & -2 &\cdot  & \cdot & \cdot  &\cdot \\
\hline
\cdot  & 0 & \cdot & \cdot  &0  & \cdot & \cdot  &\cdot  & \cdot & \cdot  &\cdot  & \cdot & -1 & 1 &\cdot  & \cdot \\ 
\cdot  &\cdot  & \cdot & \cdot  &\cdot  & \cdot & \cdot  &\cdot  & \cdot & \cdot  &\cdot  & \cdot &  0 & -1  &\cdot  & \cdot\\
\hline
\cdot  &\cdot  & \cdot & \cdot  &\cdot  & \cdot & \cdot  &-1  & \cdot & \cdot  &\cdot  & \cdot &\cdot  & \cdot&  1 & 0  \\ 
\cdot  &\cdot  & \cdot & \cdot  &\cdot  & \cdot & \cdot  &\cdot  & \cdot & \cdot  &\cdot  & \cdot &\cdot  & \cdot & -1 & 1 
\end{array}\right]
\end{align*}
(For the reader's convenience,  we have used $\cdot$ to denote a 0-entry  that occurs when $x_i$ and $x_j$ are disjoint on $F$ and geometrically unlinked, and 0 to denote a less obvious trivial linking between $x_i$ and $x_j^+$.)

%When $(n_0,n_1,n_2,n_3,n_4)=(-2,-2,2,1,-1)$,  
Column moves simplify $A+A^T$ to the following matrix:

\begin{align*}
B= \left[
\begin{array}{cccccccccccccccc}
1&0&0&0&0&0&0&0&0&0&0&0&0&0&0&0 \\
0&1&0&0&0&0&0&0&0&0&0&0&0&0&0&0 \\
0&0&1&0&0&0&0&0&0&0&0&0&0&0&0&0 \\
0&0&0&1&0&0&0&0&0&0&0&0&0&0&0&0 \\
0&0&0&0&1&0&0&0&0&0&0&0&0&0&0&0 \\
0&0&0&0&0&1&0&0&0&0&0&0&0&0&0&0 \\
0&0&0&0&0&0&1&0&0&0&0&0&0&0&0&0 \\
0&0&0&0&0&0&0&1&0&0&0&0&0&0&0&0 \\
0&0&0&0&0&0&0&0&1&0&0&0&0&0&0&0 \\
0&0&0&0&0&0&4&0&0&9&0&0&0&0&0&0 \\
0&0&0&0&0&0&0&0&0&0&1&0&0&0&0&0 \\
4&8&5&8&8&0&5&0&5& 0 & 0 & 9 & 0 & 0 & 0 & 0\\
0&0&0&0&0&0&0&0&0&0&0&0&1&0&0&0 \\
0 & 1 & 8&0&4&7&7&0&0&0&0&0&7&9&0&0\\
0&0&0&0&0&0&0&0&0&0&0&0&0&0&1&0 \\
1 & 2 & 6 &2 &2 &5 & 8 & 6 & 0 & 0 & 5 & 0 & 0 & 0 & 7 &9
\end{array}
\right]
\end{align*}
In particular,  we have established our first claim, that $H_1(\Sigma_K) \cong (\Z/ 9 \Z)^4$.  Moreover, we have that $H_1(\Sigma_K)$  is generated by $y_9, y_{11}, y_{13}, y_{15}$,  and the columns of $B$ allow us to express the other $y_i$s in terms of these generators:
\begin{align*}
\begin{array}{lllll}
y_0= 5 y_{11}+8 y_{15} && y_1= y_{11}+ 8 y_{13} +7 y_{15} && y_2= 4 y_{11}+y_{13}+3 y_{15}\\
y_3= y_{11}+ 7 y_{15} && y_4= y_{11}+5 y_{13}+7 y_{15} && y_5= 2 y_{13}+ 4 y_{15}\\
y_6= 5 y_9+4 y_{11}+2 y_{13}+y_{15} && y_7= 3 y_{15} && y_8= 4 y_{11} \\
y_{10}= 4 y_{15} && y_{12}=2 y_{13} && y_{14}= 2 y_{15}
\end{array}
\end{align*}

We now compute the homology classes of the lifts of the $\gamma_i$ and $\gamma_i'$ curves to $\Sigma_K$:
\begin{align*}
z_0 &= y_9 \\
z_1 &= y_4-y_{10}= ( y_{11}+5 y_{13}+7 y_{15})- 4 y_{15}= y_{11}+5 y_{13}+3 y_{15}\\
z_1' &= y_4=  y_{11}+5 y_{13}+7 y_{15}\\
z_2&= y_{10}-y_{11}= 8 y_{11}+ 4 y_{15}\\
z_2' &= y_4-y_5 = (y_{11}+5 y_{13}+7 y_{15})- (2 y_{13}+ 4 y_{15})=y_{11}+3 y_{13}+3 y_{15}\\
z_3 &= y_8= 4 y_{11}\\
z_3' &=y_2= 4 y_{11}+y_{13}+3 y_{15}\\
z_4 &= y_6=  5 y_9+4 y_{11}+2 y_{13}+y_{15}\\
z_4' & =y_0= 5 y_{11}+8 y_{15}
\end{align*}
Our desired conditions on $\chi$ can therefore be rewritten as follows: 
\begin{align*}
\begin{array}{llc}
(\text{C}_0) & \chi(z_0)=0: & \chi(y_9)=0 \\
(\text{C}_1) & \chi(z_1)= \pm \chi(z_1'):  & \chi(y_{15})=0 \text{ or } \chi(y_{15})= 7 \chi(y_{11})+ 8 \chi(y_{13})\\
(\text{C}_2) & \chi(z_2)= \pm \chi(z_2'): & \chi(y_{15})= 6 \chi(y_{13}) \text{ or } \chi(y_{15})= 2 \chi(y_{11})+ 3 \chi(y_{13}) \\
(\text{C}_3)& \chi(z_3)= \pm \chi(z_3'): &  6 \chi(y_{15})= \chi(y_{13}) \text{ or } 6 \chi(y_{15}) = \chi(y_{11})+ 8 \chi(y_{13}) \\
(\text{C}_4) & \chi(z_4)= \pm \chi(z_4'): & \chi(y_{13})= 2 \chi(y_9) \text{ or } \chi(y_{15})= 7 \chi(y_9)+ 4 \chi(y_{11})+ \chi(y_{13})
\end{array}
\end{align*}

Up to rescaling,  under condition $(\text{C}_0)$ any nontrivial map $\chi \colon H_1(\Sigma_K) \to \Z/ 3 \Z$ has 
\[(\chi(y_9),  \chi(y_{11}), \chi(y_{13}), \chi(y_{15}))= (0,a,b,c),\]
 where the first nonzero entry (which may or may not be $a$) can be assumed to equal $1$.  That is,  up to rescaling $(\chi(y_{11}), \chi(y_{13}), \chi(y_{15}))$ appears on the following list:
\begin{align*}
\begin{array}{c}
(0,0,1),  (0,1,0), (0,1,1),  (0,1,2), 
   (1,0,0),  (1,0,1), (1,0,2),\\
   (1,1,0) ,  (1,1,1) , (1,1,2), (1,2,0) ,  (1,2,1) , (1,2,2)
   \end{array}
\end{align*}
We now rule out all but one of these options using our conditions:
\begin{itemize}
\item Ruled out by $(\text{C}_1)$:
$(0,0,1),  (0,1,1), (1,0,2),  (1,1,1),  (1,1,2), (1,2,1)$.
\item  Ruled out by $(\text{C}_2)$:
$(0,1,2),   (1,0,1)$
\item Ruled out by $(\text{C}_3)$:
$(0,1,0), (1,2,0),  (1,2,2)$
\item Ruled out by $(\text{C}_4)$:
$(1,1,0)$. 
\end{itemize}
We have thereby established our second claim: up to rescaling there is a single nontrivial map $\chi \colon H_1(\Sigma_K) \to \Z/ 3 \Z$ satisfying $(\text{C}_0)$--$(\text{C}_4)$,  the one with $\chi(y_9)=\chi(y_{13})= \chi(y_{15})=0$ and $\chi(y_{11})=1$.  

We now wish to see if this character lifts to a surjective map $\rho$  to $\Z/ 9 \Z$.  Up to rescaling,  and assuming $(\text{C}_0)$ is satisfied,  we have that $(\rho(y_{11}), \rho(y_{13}), \rho(y_{15}))$ is of the form $(3a, 1, 3b)$ for some $a,b \in \{0,1,2\}$.  
However, all of these options are ruled out:
\begin{itemize}
\item Ruled out by $(\text{C}_1)$:
$(0,1, 3 ),  (0,1,6), (3,1,3),  (3,1,6),  (6,1,3),  (6,1,6)$
\item Ruled out by $(\text{C}_3)$:
$(0,1,0), (3,1,0), (6,1,0)$
\end{itemize}
thereby completing the proof of our third and final claim. 
\end{proof}

We will need the following version of Gilmer's Casson-Gordon signature 4-genus bound. 
\begin{theorem}[\cite{GilmerGenus}]\label{thm:Gilmer}
Let $K$ be a knot and suppose that $g_4(K) \leq g$.  Then there is a decomposition of $H_1(\Sigma_K) \cong A_1 \oplus A_2$ such that 
\begin{enumerate}
\item $A_1$ is generated by $2g$ elements
%$A_1$ has an even presentation of rank $2g$ with signature $\sigma(K)$.
\item $A_2$ has a subgroup $B$ such that  $|B|^2= |A_2|$ and
\[|\sigma(K,\chi) + \sigma(K) | \leq 4g\]
for any prime power order $\chi \colon H_1(\Sigma_K) \to \Z/q\Z$ that vanishes on $A_1 \oplus B$.
%Changed \sigma_1 \tau (K,\chi) to \sigma(K,\chi) for simplicity. 
\end{enumerate}
\end{theorem}

%and $H_1(\Sigma_2(K_{\mathbf{J}})) \cong (\Z/ 9 \Z)^4$. 
% In this case, we have that $|A_1| \leq 81$ and so $|A_2| \geq 81$ and $|B| \geq 9$.  So under the assumption that $g_4(K) \leq 1$ there exists either a surjective character $\chi \colon H_1(\Sigma_2(K_{\mathbf{J}})) \to \Z/9 \Z$ or two linearly independent characters $\chi_i \colon H_1(\Sigma_2(K_{\mathbf{J}})) \to \Z/3 \Z$  which vanish on $A_1 \oplus B$,  and hence which have ``small'' CG signature. 
We now combine Lemma~\ref{lemma:charactertracking} and Theorems~\ref{thm:Litherland} and~\ref{thm:Gilmer} to prove Proposition~\ref{prop:g4is2}. 
\begin{proof}[Proof of Proposition~\ref{prop:g4is2}]
Let $c_0:= \max\{ | \sigma(K, \chi)+ \sigma(K)|: \chi \colon H_1(\Sigma_K) \to \Z/ q \Z\}$, 
%\footnote{AM: This is sensible because there are only finitely many $\chi$.}
and let $c= \max\{c_0,2\}$.  
For $0 \leq i \leq 4$,  let $J_i$ be a knot with $|\sigma_{J_i}(\omega_9^j)|= 2^{5i+j} c$ for each $1 \leq j \leq 4$.\footnote{For example, a quick computation using the additivity of the Tristram-Levine signature function under connected sum together with the well-known formulae for the signatures of torus knots~\cite{LitherlandTorus} shows that one may take $J_i$ to be the $(2^{5i}c)$-fold connected sum of 
$3 T_{2,3} \,\# \, 3 T_{2,5} \,\# \,T_{2,7}\,\# \, 5\overline{T_{2,9}}
$.}
Let $\mathbf{J}= (J_0, J_1,J_2,J_3,J_4)$, and assume for the sake of contradiction that $g_4(K_{\mathbf{J}}) \leq 1$. 

By Theorem~\ref{thm:Gilmer},  there is a decomposition $H_1(\Sigma_{K_{\mathbf{J}}}) \cong A_1 \oplus A_2$ such that  $A_1$ is presented by a $2\times 2$ matrix and $A_2$ has a subgroup $B$ satisfying $|B|^2= |A_2|$ such that if $\chi \colon H_1(\Sigma_{K_{\mathbf{J}}}) \to \Z/ q \Z$ is a character of prime power order vanishing on $A_1 \oplus B$,  then 
\begin{align}\label{eqn:gilmerbound} |\sigma(K, \chi) + \sigma(K)| \leq 4.\end{align}

By the first claim of Lemma~\ref{lemma:charactertracking},  we have that $H_1(\Sigma_{K_\mathbf{J}}) \cong H_1(\Sigma_K) \cong (\Z/9 \Z)^4$,  and so $|A_1| \leq 81$, since the largest subgroup of $(\Z/9 \Z)^4$ presented  by a $2 \times 2$ matrix is $(\Z/9 \Z)^2$.
Therefore,  
\[|A_1 \oplus B| = |A_1| \,|B|= |A_1| \,|A_2|^{1/2}=|A_1| \left(\frac{|H_1(\Sigma_{K_{\mathbf{J}}})|}{|A_1|}\right)^{1/2}= |A_1|^{1/2}\, |H_1(\Sigma_{K_{\mathbf{J}}})|^{1/2} \leq 9^3. 
\]
In particular, there are at least $9= |H_1(\Sigma_{K_{\mathbf{J}}})|/9^3$  characters on $H_1(\Sigma_{K_{\mathbf{J}}})$ that vanish on $A_1 \oplus B$.  
Therefore,  there is either a surjective character $ H_1(\Sigma_{K_{\mathbf{J}}}) \to \Z_9$ or two linearly independent characters $H_1(\Sigma_{K_{\mathbf{J}}}) \to \Z_3$ that vanish on $A_1 \oplus B$,  and hence for which 
Equation~\ref{eqn:gilmerbound} applies. 
We say $\chi$ satisfies $(\text{C}_0)$ if $\chi(z_0)=0$ and $(\text{C}_i)$ if $\chi(z_i) = \pm \chi(z_i')$, $1 \leq i \leq 4$.  
By Lemma~\ref{lemma:charactertracking} and the previous discussion, there is a character $\chi$ that vanishes on $A_1 \oplus B$ and yet does not satisfy  $(\text{C}_j)$ for some $0 \leq j \leq 4$.  We will show that not satisfying $(\text{C}_j)$  implies that
\[ |\sigma(K_{\mathbf{J}}, \chi) + \sigma(K_{\mathbf{J}})| > 4,\]
and hence obtain our desired contradiction with Equation~\ref{eqn:gilmerbound}.

Note that $\sigma(K_{\mathbf{J}})= \sigma(K)$,  since $K_\mathbf{J}$ is obtained from $J$ by a sequence of winding number zero satellite operations. 
We now apply Theorem~\ref{thm:Litherland} to see that we have 
\begin{align}\label{eqn:CGsigformula}
 \sigma(K_{\mathbf{J}}, \chi \circ \alpha)= 
\sigma(K,  \chi) + 
2 \sigma_{J_0}(\omega_q^{\chi(z_0)})
+ 2 \sum_{i=1}^4 \left(\sigma_{J_i}(\omega_q^{\chi(z_i)})- \sigma_{J_i}(\omega_q^{\chi(z_i')})\right),
\end{align}

Let $0 \leq j \leq 4$ be maximal such that $\chi$ does not satisfy  $(\text{C}_j)$.  
If $j=0$,  we are done,  since under that assumption 
 \begin{align*}
 |\sigma(K_{\mathbf{J}}, \chi) + \sigma(K_{\mathbf{J}})| 
 %&= |\sigma_1 \tau(K,  \chi) + 2 \sigma_{J_0}(\omega_q^{\chi(z_0)})+ 2 \sum_{i=1}^4 \left(\sigma_{J_i}(\omega_q^{\chi(z_i)})- \sigma_{J_i}(\omega_q^{\chi(z_i')})\right) + \sigma(K)| 
 &= |\sigma(K,  \chi) + \sigma(K)+ 2 \sigma_{J_0}(\omega_q^{\chi(z_0)})|\\
& \geq |2 \sigma_{J_0}(\omega_q^{\chi(z_0)})|- |\sigma(K,  \chi) + \sigma(K)| 
 \geq 4c- c_0 \geq 3c  >4.  
\end{align*}

So assume now that $j \geq 1$.  Then 
\begin{align*}
(*)&=
 |\sigma(K_{\mathbf{J}}, \chi) + \sigma(K_{\mathbf{J}})| \\
 &= |\sigma(K,  \chi) + 
2 \sigma_{J_0}(\omega_q^{\chi(z_0)})
+ 2 \sum_{i=1}^4 \left(\sigma_{J_i}(\omega_q^{\chi(z_i)})- \sigma_{J_i}(\omega_q^{\chi(z_i')})\right) + \sigma(K)| \\
&\geq 2 |\sigma_{J_j}(\omega_q^{\chi(z_j)})- \sigma_{J_j}(\omega_q^{\chi(z_j')})| - | \sigma(K,  \chi)  + \sigma(K)| -
2| \sigma_{J_0}(\omega_q^{\chi(z_0)})|
-2 \sum_{i=1}^{j-1} |\sigma_{J_i}(\omega_q^{\chi(z_i)})- \sigma_{J_i}(\omega_q^{\chi(z_i')})|
\end{align*} 
Going term by term,  we have 
\begin{align*}
2|\sigma_{J_j}(\omega_q^{\chi(z_j)})- \sigma_{J_j}(\omega_q^{\chi(z_j')})| &\geq 2(2^{5j+2}- 2^{5j+1})c= 2^{5j+2}c \\
 | \sigma(K,  \chi)  + \sigma(K)|& \leq c_0 \leq c  \\
  2| \sigma_{J_0}(\omega_q^{\chi(z_0)})| &\leq 2^{5}c \\
 2|\sigma_{J_i}(\omega_q^{\chi(z_i)})- \sigma_{J_i}(\omega_q^{\chi(z_i')})|& \leq  2^{5i+5}c
 \text { for } 1 \leq i \leq j-1
\end{align*}
So 
\begin{align*}
|\sigma(K_{\mathbf{J}}, \chi) + \sigma(K_{\mathbf{J}})| 
 %& \geq  2^{5j+2}c - c- 2^{5}c- \sum_{i=1}^{j-1}  2^{5i+5}c\\
 &\geq (2^{5j+2} - 2^{5(j-1)+5}- 2^{5(j-2)+5} - \dots -2^{5}-1) c  \geq 2^{5j+1}c >4,
 \end{align*}
 and we are done. 
\end{proof}

\subsection{Arbitrarily large gaps}~\label{sec:biggaps}
We will now prove Theorem~\ref{thm:biggaps}, using the following generalization of Definition~\ref{defn:ll7strand}. 

\begin{definition}\label{defn:generalizedLL}
Given a $(6k+1)$-component string link  $\mathcal{L}$,  we define $L_{\mathcal{L}}$ as in Figure~\ref{fig:genlinkL}. 
We will only work with string links $\mathcal{L}$ that satisfy the following two properties:
\begin{enumerate}
    \item The link $L_\mathcal{L}$ shown in the Figure~\ref{fig:genlinkL} is coherently oriented and has $k+1$ components. 
    \item The two diamond markers are on the same component of $L_{\mathcal{L}}$. 
\end{enumerate}
\begin{figure}[h!]
\begin{tikzpicture}
%\draw[step=1cm,color=gray] (0,0) grid (3,5);
\node[anchor=south west, inner sep=0] at (0,0){\includegraphics[height=7cm]{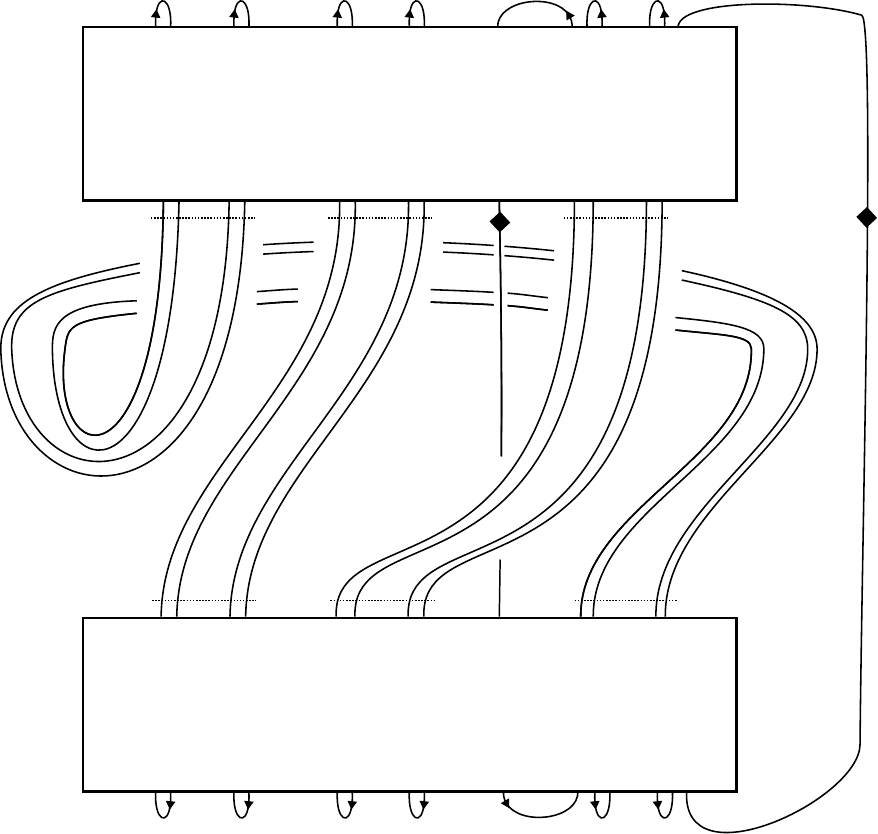}};
\node at (3.5,1){$\mathcal{L}^{-1}$};
\node at (3.4,6){$\mathcal{L}$};
\node at (1.75,5){\tiny $\dots$};
\node at (1.75,6.85){\tiny $\dots$};
\node at (1.75,.2){\tiny $\dots$};
\node at (1.1,2){\tiny $2k$};
\node at (2.6,2){\tiny $2k$};
\node at (4.6,2){\tiny $2k$};
\node at (3.23,5){\tiny $\dots$};
\node at (3.23,6.85){\tiny $\dots$};
\node at (3.23,.2){\tiny $\dots$};
\node at (5.23,5){\tiny $\dots$};
\node at (5.4,6.85){\tiny $\dots$};
\node at (5.4,.2){\tiny $\dots$};
\end{tikzpicture}
\caption{The link $L_{\mathcal{L}}$ determined by a $(6k+1)$-strand string link $\mathcal{L}$.}
\label{fig:genlinkL}
\end{figure}
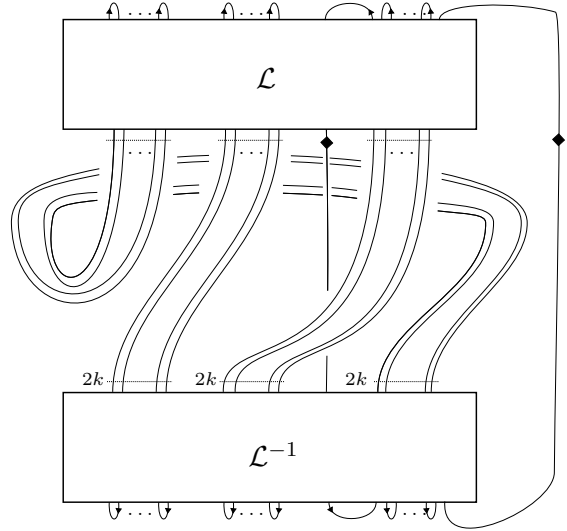
\end{definition}

The two properties of Definition~\ref{defn:generalizedLL} are exactly what is needed to mimic the argument of Proposition~\ref{prop:cobbuildingl31} and show that $L_{\mathcal{L}}$ bounds a planar surface in $\text{Spin}(L(3,1))^\circ$; we leave the details to the interested reader.

The examples for Theorem~\ref{thm:biggaps} will come from taking connected sums of knots as in Definition~\ref{defn:ll7strand}. A priori, an $m$-fold connected sum of knots as in Definition~\ref{defn:ll7strand} bounds a genus $m$ surface in $\left(\#^m\text{Spin}(L(3,1))\right)^\circ$. The following proposition shows that in fact such a knot bounds a genus $m$ surface in $\text{Spin}(L(3,1))^\circ$.

\begin{proposition}~\label{prop:connectsum}
    For $i=1,\dots, m$, let $\mathcal{L}_i$ satsify the two conditions of Definition~\ref{defn:generalizedLL} for some $k_i$, and suppose that $K_i$ is obtained from $L_{\mathcal{L}_i}$ by $k_i$ band moves. 
    Then $K= \#_{i=1}^m K_i$ bounds a surface of genus $\sum_{i=1}^m k_i$ in
    $\text{Spin}(L(3,1))^{\circ}$
\end{proposition}
Proposition~\ref{prop:connectsum} is an immediate consequence of the following lemma. 

\begin{lemma}~\label{lemm:summing}
For $i=1,2$, let $L_{\mathcal{L}_i}$ be
as in Definition~\ref{defn:generalizedLL},  where $\mathcal{L}_i$ is a string link with $6k_i+1$ components. Suppose that for $i=1,2$, the knot $K_i$ is obtained from $L_{\mathcal{L}_i}$ by $k_i$ band moves. 
Then there exists a $6k_1+6k_2+1$ component string link $\mathcal{L}$ satisfying the conditions of Definition~\ref{defn:generalizedLL} such that $K_1 \# K_2$ is obtained from $L_\mathcal{L}$ by $k_1+k_2$ band moves. 
\end{lemma}
\begin{proof}
For $\mathcal{L}_1$ and $\mathcal{L}_2$ as above, define $\mathcal{L}= \mathcal{L}_1 * \mathcal{L}_2$ to be the string link shown in Figure~\ref{fig:connsum}. Note that in Figure~\ref{fig:connsum} the box labeled $\mathcal{L}_1$ is actually the  $(6(k_1+k_2)+1)$-strand string link obtained by replacing the preferred $(4k_1+1)^{th}$ strand of $L_1$ by $(6k_2+1)$ 0-framed parallel copies, as indicated pictorally. 
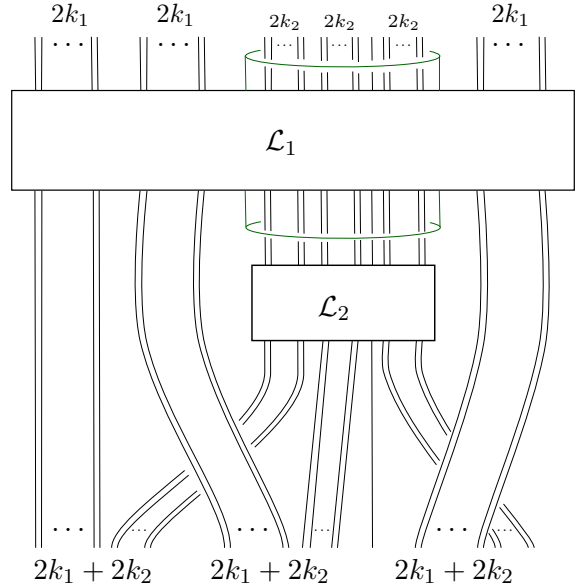
\begin{figure}[h!]
\begin{tikzpicture}
%\draw[step=1cm,color=gray] (0,0) grid (8,7);
\node[anchor=south west, inner sep=0] at (0,0){\includegraphics[height=6.8cm]{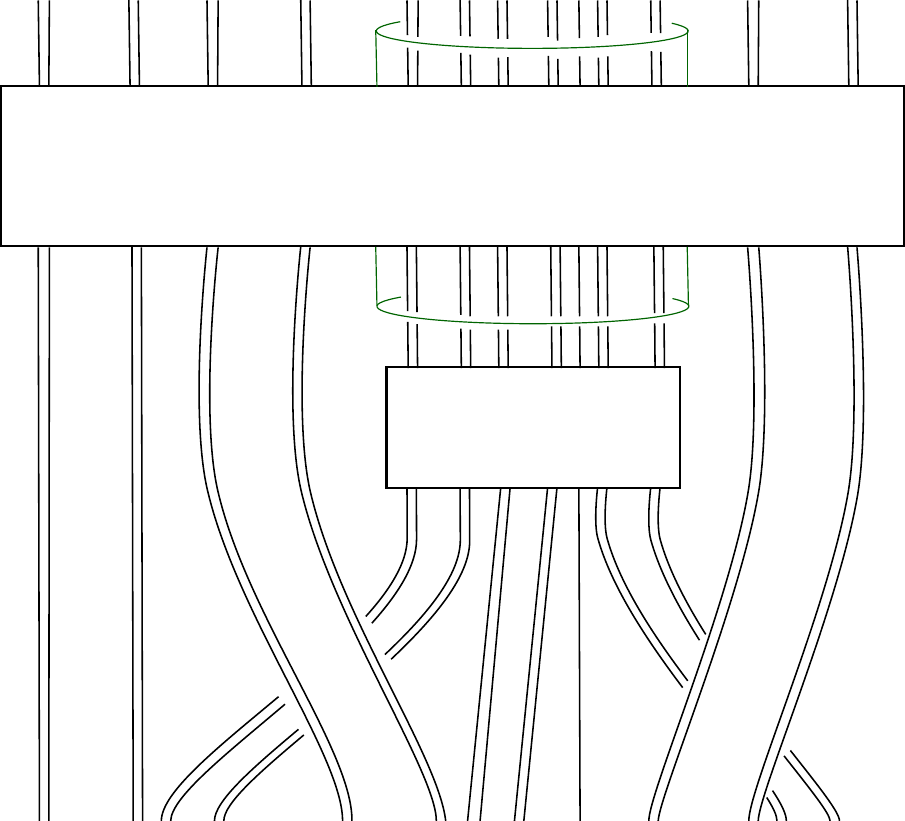}};
\node at (.8,.25) {$\dots$};
\node at (3.25,.25) {$\dots$};
\node at (5.9,.25) {$\dots$};
\node at (1.72,.25){\tiny $_{\dots}$};
\node at (4.15,.25){\tiny $_{\dots}$};
\node at (6.58,.25){\tiny $_{\dots}$};
\node at (.8,6.7) {$\dots$};
\node at (.8,7.1) {\small $2k_1$};
\node at (2.2,6.7) {$\dots$};
\node at (2.2, 7.1){\small $2k_1$};
\node at (6.65,6.7) {$\dots$};
\node at (6.65, 7.1){\small $2k_1$};
\node at (3.65,6.7){\tiny $_{\dots}$};
\node at (3.65,7){\tiny $2k_2$};
\node at (4.38,6.7){\tiny $_{\dots}$};
\node at (4.38,7){\tiny $2k_2$};
\node at (5.23,6.7){\tiny $_{\dots}$};
\node at (5.23,7){\tiny $2k_2$};
\node at (3.6,5.4){$\mathcal{L}_1$};
\node at (4.3,3.2){$\mathcal{L}_2$};
\node at (1.1, -.3){$2k_1+2k_2$};
\node at (3.45, -.3){$2k_1+2k_2$};
\node at (5.9, -.3){$2k_1+2k_2$};
\end{tikzpicture}
\caption{The $6(k_1+k_2)+1$ component string link $\mathcal{L}=\mathcal{L}_1 * \mathcal{L}_2$.}
\label{fig:connsum}
\end{figure}

One can quickly verify diagrammatically that $L_{\mathcal{L}}$ is a connected sum of $L_{\mathcal{L}_1}$ and $L_{\mathcal{L}_2}$ along the components containing the $(4k_1+1)^{th}$ strand of $\mathcal{L}_1$ and the $(4k_2+1)^{th}$ strand of $\mathcal{L}_2$. 
Additionally, the $k_1$ bands taking $L_{\mathcal{L}_1}$ to $K_1$ and the $k_2$ bands taking $L_{\mathcal{L}_2}$ to $K_2$ can be simultaneously applied to take $L_{\mathcal{L}}$ to $K_1 \# K_2$. 
\end{proof}

We are now ready to prove Theorem~\ref{thm:biggaps}. 

\begin{proof}[Proof of Theorem~\ref{thm:biggaps}]
We will take $\mathcal{K}= \#_{k=1}^{2m} K_{\mathbf{J}_k}$, where 
$\mathbf{J}_k=(J_0^k,J_1^k,J_2^k,J_3^k,J_4^k)$ are chosen to have certain Tristram-Levine signature functions. 
Since for each $k$ there is a 7-strand string link $\mathcal{L}_k$ such that $K_{\mathbf{J}_k}$ is obtained from $L_{\mathcal{L}_k}$ by a single band move, repeated application of  Proposition~\ref{lemm:summing} implies that $\mathcal{K}$ bounds a genus $2m$ surface in $\text{Spin}(L(3,1))^\circ$. 
%Note that for any choice of the $\mathbf{J}_k$ we have that $\mathcal{K}$ bounds a genus $2m$ surface in the rational homology ball that is the $2m$-fold boundary connected sum of $W$ with itself.  
So we will be done if we show that we can choose the $\mathbf{J}_k$ so that $\mathcal{K}$ has 4-genus at least $3m$.  

If $g_4(K(m)) \leq 3m-1$,  then there is an isomorphism 
$ H_1(\Sigma_{\mathcal{K}}) \cong A_1 \oplus A_2$
where $A_1$ is presented by a $(6m-2) \times (6m-2)$ square matrix and $B$ has a subgroup satisfying $|B|^2= |A_2|$  such that if $\chi$ is a character vanishing on $A_1 \oplus B$,  then 
\[ |\sigma(\mathcal{K}, \chi)+ \sigma(\mathcal{K})| \leq 4(3m-1).\]
Note that $H_1(\Sigma_{\mathcal{K}}) \cong H_1(\Sigma_K)^{2m} \cong (\Z/9\Z)^{8m}$,  and  since $A_1$ must be isomorphic to a subgroup of $(\Z/9\Z)^{8m}$ and be generated by some $(6m-2)$  elements,   its order is at most $9^{6m-2}$.  
Therefore,    there are at least
\[3^{2m+2}=\frac{3^{8m}}{3^{6m-2}}\geq \frac{|H_1(\Sigma_\mathcal{K})|^{1/2}}{|A_1|^{1/2}} = \frac{|H_1(\Sigma_\mathcal{K})|}{|A_1| \cdot\left(
\frac{|H_1(\Sigma_\mathcal{K})|}{|A_1|} \right)^{1/2}} =
\frac{|H_1(\Sigma_\mathcal{K})|}{|A_1| \cdot |A_2|^{1/2}} =
\frac{|H_1(\Sigma_\mathcal{K})|}{|A_1 \oplus B|}\]
characters vanishing on such $A_1 \oplus B$.  
 So we will be done if we can show that we can choose the $\mathbf{J}_k$ so that there are strictly fewer than $3^{2m+2}$ characters $\chi$ for which 
\begin{align}\label{eqn:gilmerboundsum}
(*)_\chi:=|\sigma(\mathcal{K}, \chi)+ \sigma(\mathcal{K})| \leq 4(3m-1).
\end{align}

We have that 
\begin{align*}
 (*)_\chi&=|\sigma(\mathcal{K}, \chi)+ \sigma(\mathcal{K})|\\
 &= 
 \left|\sum_{k=1}^{2m}\left( \sigma(K_{\mathbf{J}_k}, \chi_k) + \sigma(K)\right)\right|\\
 &= \left| 
 \sum_{k=1}^{2m}
 \left(
 \sigma(K, \chi_k)+ \sigma(K) + 2 \sigma_{J_0^k} (\omega_q^{\chi_k(z_0(k))})
+ \sum_{j=1}^4 2(\sigma_{J_i^k}(\omega_q^{\chi_k(z_i(k))}- \sigma_{J_i^k}(\omega_q^{\chi_k(z_i'(k))})
 \right)
 \right|,
\end{align*}
where $z_i(k)$ refers to the copy of $z_i$ within the $k$th copy of $H_1(\Sigma_K)$ in $H_1(\Sigma_{\mathcal{K}}) \cong H_1(\Sigma_K)^{2m}$, and similarly for $z_i'(k)$. 

As before,  we can choose our $\mathbf{J}_k$ to have sufficiently large and complex Tristram-Levine signature functions so that the only way that Equation~\ref{eqn:gilmerboundsum} is satisfied
is if for all $k$ we have $\chi_k(z_0(k))=0$ and $\chi_k(z_i(k))= \pm \chi_k(z_i'(k))$ for $1 \leq i \leq 4$. 
For example,  we can choose $J^k_i$ to be a knot with 
$\sigma_{J^k_i}(\omega_9^j)= 12mc\cdot 2^{25k+5i+j}$ for  $1\leq j \leq 4$,  where as before 
\[c= \max\{ \max\{|\sigma(K,\chi) + \sigma(K)|: \chi \colon H_1(\Sigma_K) \to \Z/q\Z\}, 2\}.\]

This completes the proof, since by Lemma~\ref{lemma:charactertracking} there are at most $3^{2m}$ characters
 \[\chi=(\chi_1,\dots, \chi_{2m}) \colon H_1(\Sigma_{\mathcal{K}})  \cong H_1(\Sigma_K)^{2m} \to \Z/q \Z
 \] such that
  each $\chi_k$ has $\chi_k(z_0(k))=0$ and 
  $\chi_k(z_i(k))= \pm \chi(z_i'(k))$ for $1 \leq i \leq 4$.   
\end{proof}
\bibliographystyle{amsalpha}
\bibliography{rationalgenus}

\end{document}